\documentclass[a4paper,12pt,draft]{article}

\usepackage{amsmath}
\usepackage{amssymb}
\usepackage{amsthm}
\usepackage{mathtools}
\usepackage{empheq}
\usepackage[final]{hyperref}
\hypersetup{
	colorlinks=true,
	citecolor=black,
	filecolor=black,
	linkcolor=blue,
	urlcolor=blue
}

\newtheorem{theorem}{Theorem}[section]
\newtheorem{lemma}[theorem]{Lemma}
\newtheorem{proposition}[theorem]{Proposition}

\newtheorem{corollary}[theorem]{Corollary}

\theoremstyle{remark}
\newtheorem{remark}[theorem]{Remark}

\def\R{\mathbb{R}}
\def\d{\partial}
\def\dif{{\mathrm d}}
\def\ut{\tilde{u}}
\def\rt{\tilde{\rho}}
\def\vt{\tilde{v}}

\def\ul{\underline}
\def\ol{\overline}

\def\hv{\hat{v}}
\def\hu{\hat{u}}
\def\mL{\mathcal{L}}
\def\mI{\mathcal{I}}
\def\mJ{\mathcal{J}}
\def\mB{\mathcal{B}}
\def\vE{\mathcal{E}}
\def\supp{\text{supp}}
\def\hb#1{{\cal B}^{#1+\sigma}}

\def\lp#1#2{|#2|_{#1}}
\def\lpp#1{\|#1\|}

\newcommand{\qu}{\quad}

\begin{document}

%% title
%%

\title{
Large-time behaviour of the spherically symmetric solution
to an outflow problem for isentropic model of compressible viscous fluid%
%\thanks{The second author's work was
%        supported in part by
%        Grant-in-Aid for Scientific Research (C)(2)
%        14540200 of the Ministry of Education, Culture, Sports,
%        Science and Technology
%          Science and Technology.}
}
\author{
Yucong {\sc Huang}${}^1$ and Shinya {\sc Nishibata}${}^2$
\\[5pt]
${}^1$ School of Mathematics\\
University of Edinburgh\\
Edinburgh EH9 3FD, UK
\\[3pt]
${}^2$ Department of Mathematical and Computing Sciences\\
Tokyo Institute of Technology\\
Tokyo 152-8552, Japan
}

\date{}
\maketitle

\vspace{-6mm}

%% abstract
%%
\begin{abstract}
We study the large time behaviour of a spherically
symmetric motion of out-flowing isentropic and compressible viscous gas. 
The fluid occupies an unbounded exterior domain in $\R^n \; (n \ge 2)$, and it flows out from an inner sphere centred at the origin of radius $r=1$.
The unique existence of a stationary solution satisfying the outflow boundary condition has been obtained in \cite{H-M}. 
The main aim of present paper is to show that this stationary solution becomes a time asymptotic state 
to the initial boundary value problem with the same boundary and spatial asymptotic conditions. 
Here, the initial data is chosen arbitrarily large if it belongs
to the suitable weighted 
Sobolev space. 
The main strategy is to approximate the unbounded exterior problem 
by solving a sequence of outflow-inflow initial boundary value problems posed in finite annular domain. 
Then the solution is obtained as a limit of these approximate solutions. 
The key argument for the stability theorem is based on the derivation of a-priori estimates in the weighted Sobolev
space, 
executed under the Lagrangian coordinate. 
The essential step of the proof is to obtain the point-wise upper and lower bound for the density. 
It is derived through employing a representation formula of the density 
with the aid of the weighted energy method.
%The H\"{o}lder regularity of the initial data is also required %for translating the results in the Lagrangian coordinate to %those in the Eulerian coordinate.
\end{abstract}
%\vspace{-3mm}
\paragraph{Keywords.}
Navier--Stokes equation, Outflow Problem, Stationary wave, Eulerian coordinate,\\
\hspace*{2.5cm}
Lagrangian coordinate, H\"{o}lder continuity.
\vspace{-6mm}
\paragraph{AMS subject classifications.}
35B35, 35B40, 76N15.

\clearpage

{\hypersetup{linkcolor=blue}
	\tableofcontents
}

\section{Introduction}
\setcounter{equation}{0}
The Navier-Stokes equation for the isentropic motion of compressible viscous gas in the Eulerian coordinate is given by
\begin{subequations}
\label{ns}
\begin{gather}
\rho_t + \nabla \cdot (\rho u) = 0, \label{ns1} \\
\rho \{ u_t + (u \cdot \nabla) u \} =
\mu_1 \Delta u +(\mu_1 + \mu_2) \nabla (\nabla \cdot u)
- \nabla P(\rho).
\label{ns2}
\end{gather}
\end{subequations}
We study the asymptotic behaviour of a solution $(\rho, u)$ to (\ref{ns}) in an unbounded exterior domain $\Omega \vcentcolon= \{ z \in \R^n \; ; \; |z| > 1 \}$, where $n$ is the space dimension larger than or equal to $2$. Here $\rho>0$ is the mass density; $u=(u_1,\dots,u_n)$ is the velocity of gas; $P(\rho) = K\rho^\gamma \; (K>0, \gamma \ge 1)$ is the pressure
with the adiabatic exponent $\gamma$; $\mu_1$ and $\mu_2$ are constant called viscosity-coefficients satisfying $\mu_1>0$ and $2\mu_1 + n \mu_2 \ge 0$. In the equations (\ref{ns}), we use notations:
$\nabla \cdot u \vcentcolon= \sum_{i=1}^n \d_i u_i$,
$(u \cdot \nabla) u \vcentcolon= ((u \cdot \nabla) u_1, \dots, (u \cdot \nabla) u_n)$,
$(u \cdot \nabla) u_j \vcentcolon= \sum_{i=1}^n u_i \d_i u_j$,
$\Delta u \vcentcolon= (\Delta u_1, \dots, \Delta u_n)$,
$\Delta u_j \vcentcolon= \sum_{i=1}^n \d_i^2 u_j$,
$\nabla P \vcentcolon= ( \d_1 P, \dots, \d_n P )$
and $\d_i \vcentcolon= \d_{z_i}$.

It is assumed that the initial data is also spherically symmetric. Namely, for  $r \vcentcolon= |z|$
\begin{equation*}
\rho_0(z) = \hat{\rho}_0(r), \quad u_0(z) = \dfrac{z}{r} \hat{u}_0(r).
\end{equation*}
Under these assumptions, it is shown in \cite{itaya85} 
that the solution $(\rho, u)$ is spherically symmetric. Here, the spherically symmetric solution is a solution to (\ref{ns}) in the form of
\begin{equation}
\rho(z,t) = \hat{\rho}(r,t), \quad u(z,t) = \frac{z}{r} \hat{u}(r,t).
\label{sym}
\end{equation}
Substituting (\ref{sym}) in (\ref{ns}), we reduce the system (\ref{ns}) to that of the equations for $(\hat{\rho}, \hat{u})(r,t)$, which is (\ref{nse}) below. Here and hereafter, we omit the hat ``\,$\hat{ \ }$\,'' to express  spherically symmetric functions without confusion.
Hence the spherically symmetric solution $(\rho, u)$
satisfies
the system of equations
\begin{subequations}
\label{nse}
\begin{gather}
\rho_t + \frac{(r^{n-1} \rho u)_r}{r^{n-1}} = 0, \label{nse1} \\
\rho (u_t + uu_r) = \mu \Big( \frac{(r^{n-1} u)_r}{r^{n-1}}
\Big)_r
- P(\rho)_r, \label{nse2}
\end{gather}
\end{subequations}
where $\mu := 2\mu_1 + \mu_2 $ is a positive constant,
standing for viscosity constant.
 The initial data
to (\ref{nse})
is prescribed by
\begin{subequations}
\label{iac}
\begin{gather}
\rho(r,0) = \rho_0(r)>0, \quad u(r,0) = u_0(r), \label{ic}\\
\lim_{r \to \infty} (\rho_0(r), u_0(r)) = (\rho_+, u_+), \quad 0<\rho_+ <\infty, \label{asymp}
\end{gather}
\end{subequations}
where $\rho_+$, $u_+$ are constants.
In the present paper, 
we study
the problem where fluid flows out, at a constant velocity,
through the boundary of
inner sphere, centred at the origin of radius $r=1$. The corresponding boundary condition is given by
\begin{equation}\label{ub}
 u(1,t) = u_b, \ \ \text{where \ $u_b<0$}. 
\end{equation}
Moreover, it is assumed that the initial data (\ref{ic}) is compatible with the boundary data (\ref{ub}):
\begin{subequations}\label{compa}
\begin{gather} 
u_0(1) = u_b, \label{compa1} \\ 
\Big\{ \rho_0 u_0 (u_0)_r + \mu \Big( \frac{(r^{n-1}u_0)_r}{r^{n-1}}\Big)_r - P(\rho_0)_r \Big\} \bigg|_{r=1} = 0. \label{compa2}
\end{gather}
\end{subequations}
Since the characteristic speed of (\ref{nse1}) is negative
on the boundary due to (\ref{ub}), no boundary condition for the density $\rho$ around $r=1$ is necessary for the well-posedness of the initial boundary value problem (\ref{nse}), (\ref{iac}) and (\ref{ub}).

\paragraph{\bf Notations.} Before stating the main result, we list the notations which will be frequently used in the present paper:
\begin{itemize}
\item $c$ and $C$ denote generic positive constants, which depend only on the quantities $\rho_{+}$, $\mu$, $\gamma$, $K$, and $n$. 
\item For a non-negative integer $k \ge 0$, $H^k(\Omega)$ denotes the $k$-th order Sobolev space over $\Omega$ in the $L^2$ sense with the norm
\begin{equation*}
\| f \|_k\vcentcolon= \|f\|_{H^{k}(\Omega)} = \bigg( \sum_{|\alpha|=0}^{k} \int_{\Omega}\!\!  |\d_x^{\alpha} f (x)|^2 \, \dif x \bigg)^{\frac{1}{2}}.
\end{equation*}
We note also that $H^0 = L^2$ and denote $\lpp{\cdot}  \vcentcolon= \| \cdot \|_0$.
\item For a non-negative integer $k\ge 0$, $\mB^k(\Omega)$ denotes the space of all functions $f$ which, together with all their partial derivatives $\d_x^{i}f$ of orders $|i|\le k$, are continuous and bounded on $\Omega$. It is endowed with the norm:
\begin{equation*}
|u|_{k}^{\Omega}
\vcentcolon= \sum_{|i|=0}^k \sup_{x \in \Omega} |\d_x^i u(x)|.
\end{equation*}
Moreover, for $\alpha \in (0,1)$, $\mB^\alpha (\Omega)$ denotes the
space of bounded functions over $\Omega$ which have the
uniform H\"older continuity with exponent $\alpha$.
For an integer $k\ge 0$,  $\mB^{k+\alpha}(\Omega)$ denotes the space of the functions satisfying %
$\d_x^{i} u \in \mB^\alpha (\Omega)$ for all integer $i \in [0, k]$. It is endowed with the norm:
$|\cdot|_{k+\alpha}$ is its norm defined by
\begin{equation*}
 |u|_{k+\alpha}^{\Omega}
\vcentcolon= \sum_{|i|=0}^k \sup_{x \in \Omega} |\d_x^i u(x)|
+ \sup_{\begin{smallmatrix} x,x' \in \Omega \\ x \neq x' \end{smallmatrix}}
  \frac{|\d_x^k u(x) - \d_x^k u(x')|}{|x - x'|^\alpha}.
\end{equation*}
\item For a domain $Q_T \subseteq [0,\infty) \times [0, T]$, $\mB^{\alpha, \beta} (Q_T)$ denotes the space of the uniform H\"{o}lder continuous functions with the H\"{o}lder exponents $\alpha$ and $\beta$ with respect to $x$ and $t$, respectively. For integers $k$ and $l$, $\mB^{k+\alpha, l+\beta}(Q_T)$ denotes the space of the functions satisfying $\d_x^i u, \d_t^j u \in \mB^{\alpha, \beta}(Q_T)$ for all integers $i \in [0,k]$ and $j \in [0,l]$. $| \cdot |_{k+\alpha,l+\beta}^{Q_T}$ is its norm.
\end{itemize}

\paragraph{Stationary solutions and main result.} The initial boundary value problem (\ref{nse})-(\ref{ub}) is formulated to study the behaviour of compressible viscous gas flowing out
through a boundary. We show that the time asymptotic state of solution to (\ref{nse})-(\ref{ub}) is the stationary solution, which is a solution to (\ref{nse}) independent of time
variable
$t$, satisfying the same boundary conditions (\ref{asymp}) and (\ref{ub}). Hence the stationary solution $(\rt, \ut)(r)$ satisfies the system of equations
\begin{subequations}
\label{st}
\begin{gather}
\frac{1}{r^{n-1}}(r^{n-1} \rt \ut)_r = 0, \label{st1} \\
\rt \ut \ut_r = \mu \Big( \frac{(r^{n-1} \ut)_r}{r^{n-1}} \Big)_r
- P(\rt)_r, \label{st2}
\end{gather}
\end{subequations}
for $r\in[1,\infty)$ and the 
boundary and the spatial asymptotic conditions
\begin{equation}\label{stbdry}
 \ut(1) = u_b, \quad
 \lim_{r \to \infty} (\rt(r), \ut(r)) = (\rho_+, u_+).
\end{equation}
Multiplying $r^{n-1}$ on (\ref{st1}) and integrating the resultant equality over $(1,r)$, we obtain
\begin{equation}\label{stu}
\ut(r) = u_b \rt(1) \dfrac{r^{1-n}}{\rt(r)} \quad \text{ for $r\ge1$.}
\end{equation}
Since $n\ge2$ and $\rho_+>0$, if a solution $(\rt,\ut)(r)$ exists and $\rt(1)<\infty$, then it is necessary
that $\ut(r)\to 0$ as $r\to\infty$. Thus we impose the far-field condition
\begin{equation}\label{u+}
	u_+ = 0.
\end{equation}
Substituting (\ref{stu}) in (\ref{st2}) yields the differential equation
\begin{equation}\label{stODE}
\dfrac{\dif^2\rt}{\dif r^2} + a_1(\rt,r) \Big\vert \dfrac{\dif \rt}{\dif r}\Big\vert^2 + a_2(\rt,r)\dfrac{\dif \rt}{\dif r} + a_3(\rt,r)= 0, \quad \lim\limits_{r\to\infty} \rt(r) = \rho_+.
\end{equation}
where $a_1$, $a_2$, $a_3$ are given by
\begin{gather*}
a_1(\tilde{\rho},r):= -\dfrac{2}{\rt(r)}, \qquad a_3(\rt,r):= - \dfrac{(n-1)u_b}{\mu} \dfrac{\rt(1)\rt(r)}{r^n},\\
a_2(\tilde{\rho},r) := \dfrac{r^{n-1}}{\mu u_b } \dfrac{\rt(r)P^{\prime}(\rt(r))}{\rt(1)} - \dfrac{u_b}{\mu}\dfrac{\rt(1)}{r^{n-1}} - \dfrac{n-1}{r}.
\end{gather*}
For $P(\rho)=K\rho^{\gamma}$, 
the 
equation (\ref{stODE}) is solved 
by 
using 
an iteration scheme in \cite{h-m21}. The paper \cite{h-n-s23}
obtains the decay rate of the stationary solution $(\rt,\ut)$
as $r \to \infty$.
These results are summarized as follows.
\begin{lemma}\label{lemma:st}
There exists $\delta=\delta(\rho_+,\mu,\gamma,K,n)>0$ such that if $| u_b | \le \delta$, then there exists a unique solution $(\rt,\ut)\in \mB^2[1,\infty)$ to the problem (\ref{st})--(\ref{u+}). In addition, $\rt(r)$ and $\ut(r)$ are  strictly monotone increasing and decreasing, respectively. Moreover, there exists constant $C=C(\rho_+,\mu,\gamma,K,n)>0$ such that for $r\ge1$,
\begin{subequations}\label{stProp}
\begin{gather}
r^{n-1}\rt(r)\ut(r)= \rt(1) u_b , \qu C^{-1} |u_b| r^{-n} \le \ut_r(r) \le C |u_b| r^{-n},\label{stProp1}\\ 
C^{-1} |u_b|^2 r^{-2n+2} \le  \rho_+ - \rt(r)  \le C |u_b|^2 r^{-2n+2},  \label{stProp2}\\
C^{-1} |u_b|^2 r^{-2n+1} \le \rt_r(r) \le C |u_b|^2 r^{-2n+1}, \label{stProp3}\\
%C^{-1} |u_b| r^{-2n} \le |\ut_{rr}(r)| \le C |u_b|r^{-2n}, \qu
C^{-1} |u_b|^2 r^{-2n} \le |\rt_{rr}(r)| \le C |u_b|^2 r^{-2n}.\label{stProp4}
\end{gather}
\end{subequations}
\end{lemma}
The paper \cite{h-n-s23}
also 
shows the asymptotic stability of the stationary solution
$(\rt,\ut)$
for small initial perturbation of $(\rt,\ut)$
by utilizing the decay rate \eqref{stProp3}.
The main result of the present
paper is to obtain 
the time asymptotic stability of the stationary solution $(\rt,\ut)$ for the large initial perturbation in the suitable Sobolev space. This 
result
is stated as follows.
\begin{theorem}
\label{thm:main}
Assume $n \ge 2$ and $1\le \gamma\le 2$. Let $\sigma$ be an arbitrary
constant satisfying $0 < \sigma < 1$.
Suppose the initial data (\ref{iac}) belongs to the function space
\begin{subequations}
\begin{gather}
 r^{\frac{n-1}{2}} (\rho_0 - \rt), \  r^{\frac{n-1}{2}} (u_0-\ut), \
 r^{\frac{n-1}{2}} (\rho_0 - \rt)_r, \ r^{\frac{n-1}{2}} (u_{0}-\ut)_r
 \in L^2(1, \infty), \label{Hea}
\\
 \rho_0 \in \hb{1}[1,\infty), \
 u_0\in  \hb{2} [1,\infty)
%\quad \text{for a certain} \ \sigma \in (0,1),
\label{He}
\end{gather}
\end{subequations}
and satisfies (\ref{u+}) and the compatibility condition (\ref{compa}). Then there exists $\delta_0=\delta_0(\rho_0,u_0,\rho_+,\mu,\gamma,K,n)>0$ such that, for $|u_b|\le \delta_0$, the initial boundary value
problem (\ref{nse})--(\ref{ub}) has a unique solution $(\rho, u)$ satisfying
\begin{subequations}
\label{msp}
\begin{gather}
r^{\frac{n-1}{2}} (\rho - \rt), \
r^{\frac{n-1}{2}} (u-\ut), \
r^{\frac{n-1}{2}} (\rho - \rt)_r, \
r^{\frac{n-1}{2}} (u-\ut)_r
\in C ([0, T] ; L^2(1, \infty)), 
\label{msp1} \\
\rho \in \mB^{1+\sigma,1+\sigma/2}
([1,\infty) \times [0,T]), \
u \in \mB^{2+\sigma, 1+\sigma/2}
([1,\infty) \times [0,T])
\label{msp2}
\end{gather}
\end{subequations}
for an arbitrary $T>0$. Furthermore, $(\rho, u)$ converges to the corresponding stationary solution $(\rt,\ut)$ 
given by Lemma \ref{lemma:st} as time $t$ tends to infinity. Precisely, it holds that
\begin{equation}\label{tcvg}
 \lim_{t \to \infty} \sup_{r \in [1,\infty)}
 | (\rho(r,t)-\rt(r), u(r,t)- \ut(r))  | = 0.
%\label{asym}
\end{equation}
\end{theorem}
Notice that any smallness assumptions on
the initial data is not necessary in the above stability theorem.
Owing to \cite{itaya71,Tani}, the H\"older continuity of the initial data (\ref{He}) ensures the unique 
existence of the 
H\"older continuous 
solution
locally in time. This regularity also
validates the transformation between
the Eulerian and the Lagrangian 
coordinates (see (\ref{B})--(\ref{R}) below). We then derive 
the
a-priori Sobolev estimates in the Lagrangian domain, 
and using these estimates the H\"older
continuity of $(\rho,u)$ is
attained by the Schauder theory for parabolic equations.
By these procedures,
we show the existence 
of the solution
globally in time.  Namely, 
the solution of is extended to
an arbitrary time $T>0$ 
owing to
 the continuity argument. 
 For the derivation of H\"older continuity, readers
 are refered
 to \cite{k-n-z03,NNY,NN}. 
 See also \cite{friedman64,lady} for 
 a survey on the general theory of the Schauder estimates
 for parabolic equations. 

%The remainder of the present paper is following.
%After summarizing the related results, we start the
%detailed discussions on the proof of Theorem \ref{thm1}.
%In Section 2, we derive the a priori estimate of the solution
%by employing the standard energy method.

\paragraph{Related results.}
The well-posedness of compressible
Navier-Stokes equations is an important subject in the development of mathematical physics.
Particularly, the problem 
with an out-flowing boundary condition has gained a lot of interests in the recent years. Here, we state several previous results, which are relevant to the present paper. 

First of all, for a comprehensive survey of the mathematical theory of compressible viscous flow, we refer readers to the book \cite{kaz90} by S.~N. Antontsev, A.~V. Kazhikhov, and V.~N. Monakhov. The first notable research on the large time stability of solution is obtained by A.~Matsumura and T.~Nishida in \cite{m-n83}, where they study
the equations for
the heat-conductive compressible flow for
a general $3$-dimensional exterior domain with no-flow boundary conditions (Dirichlet boundary condition for the velocity, and Neumann boundary condition for the temperature). Under smallness assumptions on the initial data and exteral forces, they proved the global-in-time stability of the stationary solution.

When the equation is spherically symmetric with no-flow boundary condition, a pioneering work has been done by N.~Itaya \cite{itaya85}, which establishes the global-in-time existence and uniqueness of a classical solution on a bounded annulus domain, without smallness assumption on the initial data. Later, T.~Nagasawa studied the asymptotic state for the same
problem in \cite{nagasawa}. The paper \cite{itaya85} has spurred a sequence of developments on the topic of spherically symmetric solution. For example, A.~Matsumura in \cite{matsu92} constructed a spherically symmetric classical solution to the isothermal model with external forces on a bounded annular domain. Moreover, he also showed the convergence to the stationary solution as time tends to infinity, with an exponential convergence rate. Subsequently, the result of \cite{matsu92} has been extended to the isentropic and heat-conductive models by K.~Higuchi in \cite{higuchi92}. The well-posedness of spherically symmetric solution in an unbounded exterior domain was first obtained by S.~Jiang \cite{jiang96}, where the global-in-time existence and uniqueness of a classical solution is shown. In addition, a partial result on the time asymptotic
stability 
is proved in \cite{jiang96} where, for $n=3$, $\|u(\cdot,t)\|_{L^{2j}}\to 0$ as $t\to\infty$ with any fixed integer $j\ge 2$. Later, this restriction on the long time stability was fully resolved by T.~Nakamura and S.~Nishibata in \cite{n-n08}, where a complete stability theorem was obtained for both $n=2$ and $3$ with large initial data. We also refer to the paper \cite{NNY} by T.~Nakamura, S.~Nishibata, and S.~Yanagi, where they 
show
the time asymptotic stability of the spherically symmetric solution for the isentropic flow with large initial data and external forces. The present paper borrows several ideas from
these researches.

For the general outflow or inflow problem, where the velocity is not zero at the boundary, the stationary solution becomes non-trivial, which leads to a variety of significantly more interesting time-asymptotic behaviours for the solutions. A.~Matsumura and K.~Nishihara in \cite{m-n01} started the first investigation of this problem for the isentropic model posed on the one dimensional half-space domain. Several kinds of boundary conditions were studied in \cite{m-n01}, which includes inflow, outflow, and no flow boundary conditions. They formulated conjectures on the classification of asymptotic behaviours of the solutions in different cases subject to the relation between the boundary data and the spatial asymptotic data. Then the stability theorems for some cases were established by employing the Lagrangian mass coordinate. Following this work, S.~Kawashima, S.~Nishibata and P.~Zhu in \cite{k-n-z03} further examined the outflow problem in half line. They obtain
a set of a priori estimates directly in Eulerian coordinate, which leads to the long time stability of solutions with small initial perturbation with respect to the stationary solution. 
The convergence rate towards the stationary
solution is obtained
by T.~Nakamura, S.~Nishibata and T.~Yuge in \cite{[NNY]}.

For the non-isentropic inflow problem in the half line, T.~Nakamura and S.~Nishibata in \cite{NN} demonstrate
the time-asymptotic stability of stationary solutions under a small initial perturbation, for both the subsonic and transonic cases. 

The research on the outflow and inflow problems for the spherically symmetric solution in an unbounded exterior domain has been relatively new. Under the assumption that speed at boundary, $|u_b|$ is sufficiently small, I.~Hashimoto and A.~Matsumura in \cite{h-m21} employed the iteration method to obtain the existence and uniqueness of spherically symmetric stationary solutions for both inflow and outflow problems in an exterior domain. More recently, I.~Hashimoto, S.~Nishibata, S.~Sugizaki in \cite{h-n-s23} showed the stability of corresponding stationary solution under a small spherically symmetric perturbation. The present paper shows the stability, without the smallness assumption
on the initial data, for the outflow problem.

\paragraph{Outline of the paper.}
The remainder of this paper is organized as follows. In Section \ref{sec:approx}, we construct a set of approximation problems posed in a finite annular domain, with outer radius parametrised by $m\in\mathbb{N}$. These problems are then reformulated from the original Eulerian coordinate to Lagrangian coordinate in Section \ref{sec:RL}. The main purpose of this is to obtain the a priori uniform point-wise bounds on density, and this is derived in Section \ref{sec:apriori}. Consequently in Section \ref{sec:higher}, the uniform estimates on the higher derivatives of approximate solutions is obtained using the point-wise bound of density. We remark that the $H^1$ estimate of $\rho$ has also utilised the special structure of the equations under 
Lagrangian coordinate. Then with the help of these uniform a priori estimates, we obtain in Section \ref{sec:limit}, solution to the original problem in unbounded domain as the limit of approximate solutions via the compactness argument. Finally, making use of the estimates obtained previously, we also show the asymptotic stability of the stationary solution in Section \ref{sec:tAsymp}.

In these discussions on the solution in the Lagrangian coordinate,
the H\"{o}lder regularity of the solution is not necessary.
It is specifically required in the
time local existence of the solution and 
translating the results in this coordinate
to those in the Eulerian coordinate. 
This is proved in Section \ref{sec:higher} by utilizing the Schauder theory for the parabolic equation.

\section{Construction of approximation solutions}\label{sec:approx}\setcounter{equation}{0}
Theorem \ref{thm:main} is proved by combining the
local existence and the a priori estimate. In order to prove the local existence to the problem
\eqref{nse}--\eqref{ub}, we solve the approximate problems of \eqref{nse} in bounded annular domains $(r,t)\in[1,m]\times[0,T]$, which are parametrised by $m\in\mathbb{N}$. This procedure is necessary since several
coefficients in \eqref{nse} are unbounded over $r \in [1, \infty)$.

For each $m\in\mathbb{N}$, we seek solutions $(\rho_m,u_m)(r,t)$ to the following approximate problem in a bounded annular region $r\in[1,m]$:
\begin{subequations}\label{approx}
\begin{align}%[left = {\empheqlbrace \,}, right = {}]{align}
&\text{$(\rho_m,u_m)(r,t)$ solves \eqref{nse}} &&  \text{in $(r,t)\in[1,m]\times[0,\infty)$,} \label{approx1}\\
&\begin{aligned}
	& u_m(1,t)=u_b,\\
	& u_m(m,t)=\ut(m), \ \ \rho_m(m,t)=\rt(m)
\end{aligned} && \text{for $t\in[0,\infty)$,}\label{approx2}\\
&(\rho_m,u_m)(r,0)= (\rho_m^0,u_m^0)(r) && \text{for $r\in[1,m]$.}\label{approx3}
\end{align}
\end{subequations}
Here, $(\rho_m^0,u_m^0)(r)$ is the modified initial data constructed from $(\rho_0,u_0)(r)$ as follows: let $\varphi_m(x) \in \mB^3[1,\infty)$ be the cut-off functions satisfying
\begin{gather}
	\varphi_m(r) \vcentcolon=
	\left\{
	\begin{array}{ll}
		1,     & \text{for} \ \ 0 \le r \le m/2, \\
		0,     & \text{for} \ \ m \le r,
	\end{array} \right. \label{cut1}\\
	0 \le \varphi_m(r) \le 1, \quad| \partial_r^i \varphi_m(r) | \le C m^{-i} \quad (i = 1,2,3) \quad \text{for} \quad r\in[1,\infty).
	\nonumber
\end{gather}
The approximate initial data $(\rho_{m}^0, u_{m}^0)$ is 
derived,
from \eqref{iac}, by using $\varphi_m(x)$ as
\begin{equation}\label{inim}
\rho_{m}^0(r) \vcentcolon= (\rho_0(r) - \rt(r)) \varphi_m(r) + \rt(r), \quad
u_{m}^0(r) \vcentcolon= (u_0(r) - \ut(r)) \varphi_m(r) + \ut(r)
\end{equation}
for $r\in[1,\infty)$. We define the weighted $H^1$ norm by
\begin{equation*}
	\|f\|_{1,r,m}
	\vcentcolon= \Big( \int_{1}^{m} \{ |f|^2 + |f_r|^2 \}(r) r^{n-1} \dif r \Big)^{\frac{1}{2}} ,\ \ E_m^0 \vcentcolon= \| (\rho_{m}^0 - \rt , u_m^0 - \ut ) \|_{1,r,m}^2.
\end{equation*}
Owing to the above construction, 
it is easy to show the following proposition.
\begin{proposition}\label{prop:icm}
There exists a constant $C_0>0$ independent of $m\in\mathbb{N}$ such that
\begin{align*}
&C_0^{-1} \le \rho_m^0(r) \le C_0 \quad \text{ for } \ \ r\in[1,\infty),  \\
&\rho_m^0(m)=\rt(m), \quad u_m^0(m)=\ut(m), \quad (\rho_m^0,u_m^0) \text{ satisfies \eqref{compa}}, \\ 
&\sup_{m\in\mathbb{N}}\| (\rho_m^0-\rt,u_m^0-\ut) \|_{1,r,m}^2 \le C_0 E_0, \ \ \text{where} \ \ E_0\vcentcolon=\|(\rho_0-\rt,u_0-\ut)\|_{1,r,\infty}^2.
%\sup_{m\in\mathbb{N}}\int_{1}^{m} \big\{ |\rho_m^0-\rt|^2 + | u_m^0 - \ut |^2 \big\}(r) r^{n-1} \dif r \le C_0,\\
%\sup_{m\in\mathbb{N}}\int_{1}^{m} \big\{ | \partial_r (\rho_m^0-\rt) |^2 + | \partial_r ( u_m^0 - \ut ) |^2 \big\}(r) r^{n-1} \dif r \le C_0.
\end{align*}
In addition, the convergence holds as $m\to\infty$:
\begin{equation*}
	(\rho_m^0,u_m^0)(r)\to (\rho_0,u_0)(r) \qu \text{for all } \ r\in[1,\infty), \qu E_m^0 \to E_0.
\end{equation*} 
\end{proposition}

The local existence of the solution to the problem \eqref{approx} in bounded domain is proved by the standard iteration method and commutator arguments. See \cite{itaya71,m-n83} for example. For $\bar{d} > \ul{d} > 0$, $D>0$ and positive integer $m$,
we define the function space as
\begin{equation*}
	\begin{aligned}
		X_{\ul{d},\bar{d},D}^m(0,T)
		 \vcentcolon= \big\{ (\rho, u) \; | \; & (\rho - \rt, u-\ut) \in C^0([0,T]; H^1(1,m)), \ \ul{d} \le \rho(x,t) \le \bar{d} \\
		&u-\ut \in L^2(0,T; H^2(1,m)), \ \| (\rho-\rt, u-\ut)(t)\|_{1,r,m} \le D  \big\},	
	\end{aligned}
\end{equation*}

\begin{lemma}\label{lemma:local}
If the initial data satisfies $ E_m^0  \le D_0$ and $\displaystyle \ul{d}_0 \le \rho_{m}^0(r) \le \bar{d}_0$ for certain constants $0<\underline{d}_0<\bar{d}_0$ and $D_0>0$ independent of $m\in\\mathbb{N}$, then there exists a positive time $T=T(\underline{d}_0,\bar{d}_0,D_0)>0$ such that the problem \eqref{approx} has a unique solution $(\rho_m, u_m)$ in the space $X_{\ul{d}_0/2,2\bar{d}_0,2D_0}^m(0,T)$.
\end{lemma}

\begin{remark}\label{rem:exT}
The existence time $T_m$ satisfies that
Since $\displaystyle \rho_{m}^0(r) > \min \{ \rho_0(r), \rt(r) \}$ and $E_m^0 \le C E_0$, the time of existence $T=T(\ul{d}_0,\bar{d}_0,D_0)$ depends only on $\displaystyle \inf_{x \in (0,\infty)} v_0(x)$ and $E_0$.
\end{remark}

\section{Reformulation in Lagrangian coordinate}\label{sec:RL}\setcounter{equation}{0}
The main strategy for the proof of Theorem~\ref{thm:main} is to employ the energy method and derive a priori estimates which are uniform over $m\in\mathbb{N}$ and $T>0$. For this purpose, it is necessary to adopt the Lagrangian mass coordinate rather than the Eulerian coordinate, since the uniform point-wise bound of $\rho$, as well as the $H^1$ estimate of $\rho$
is only  derived under Lagrangian formulation.
By Lemma \ref{lemma:local} and Remark \ref{rem:exT}, let $(\rho_m,u_m)(r,t)$ be the local-in-time approximate solution, and let $T>0$ be the maximum time of existence. Until the end of Section \ref{sec:higher}, we will denote for simplicity: $$\rho\equiv\rho_m, \qu u\equiv u_m.$$ In order to determine the coordinate transformation law for the outflow problem, we define the curve $B:[0,T]\to [0,\infty)$ as
\begin{equation}\label{B}
B(t)\vcentcolon= -u_b \int_{0}^{t} \rho(1,s) \dif s = |u_b| \int_{0}^{t} \rho(1,s) \dif s. 
\end{equation}
Physically, $B(t)$ indicates the total amount of mass flown out from the boundary $r=1$, after time $t>0$ has passed. Using this, we also define $M\vcentcolon[0,T]\to[0,\infty)$ as
\begin{equation}\label{Mdef}
	M(t)\vcentcolon= B(t) + \int_{1}^{m} \rho (r,t) r^{n-1} \dif r.
\end{equation}
By \eqref{approx1} and \eqref{approx2}, we obtain
\begin{equation*}
	M^{\prime}(t) = -u_b \rho(1,t) - \int_{1}^{m} (r^{n-1}\rho u)_r(r,t) \dif r = -m^{n-1} \rt(m) \ut(m). 
\end{equation*}
Using \eqref{stu}, it follows that $m^{n-1}\rt(m) \ut(m)= \rt(1) u_b$, hence $M^{\prime}(t)=-\rt(1)u_b=\rt(1)|u_b|$. Therefore $M(t)$ is a solution to the ordinary differential equation:
\begin{equation}\label{MODE}
	M^{\prime}(t)=\rt(1)|u_b| \qu \text{for $t\in[0,T]$}, \qu M(0)=M_0\vcentcolon=\int_{1}^{m} \rho_m^0(r)r^{n-1} \dif r>0. 
\end{equation}
By the construction \eqref{cut1}, we
verify that $M_0=M_0(m)\to\infty$ as $m\to\infty$. Since the solution to \eqref{MODE} is unique, we obtain
\begin{equation}\label{M}
	M(t) = M_0 + \rt(1) |u_b| t \ \ \text{for all $t\in[0,T]$.}
\end{equation}
Next, we consider the map $(r,t)\mapsto X(r,t)\vcentcolon [1,m]\times[0,T]\to \R$ defined by
\begin{equation}\label{RInv}
	X(r,t) \vcentcolon= B(t) + \int_{1}^{r} \rho(y,t) y^{n-1} \dif y.
\end{equation}
Since $\ul{d} \le \rho(r,t)\le \bar{d}$ in $(r,t)\in[1,m]\times[0,T]$ for some $0<\ul{d}<\bar{d}<\infty$, one can employs Inverse Function theorem to obtain that, for each $t\in[0,T]$ there exists a unique map $x\mapsto R(x,t) \vcentcolon [B(t),M(t)]\to [1,m]$ such that $X(R(x,t),t)=x$. In other words,
\begin{subequations}\label{Rexist}
\begin{gather}
	x= B(t) + \int_{1}^{R(x,t)} \rho(r,t) r^{n-1} \dif r, \quad \text{for each $B(t) \le  x\le M(t)$.}\label{Rmass}\\
	\text{Thus } \ \ R(B(t),t)=1 \ \text{ and } \ R(M(t),t)=m \ \text{ for each } \ t\in[0,T].\label{Rbdry}
\end{gather} 
\end{subequations}
Using the above construction, we set the Lagrangian space-time domain $\mL(T)$, and the Lagrangian snapshot domain $S(t)$ at time $t\in[0,T]$ to be:
\begin{subequations}\label{Ldomain}
\begin{align}
\mL(T)\vcentcolon=& \{ (x,t)\in\R\times[0,T] \;|\; B(t) \le x \le M(t) \}, \label{lsp} \\
S(t) \vcentcolon=& \{ x\in\R \;|\; B(t) \le x \le M(t) \} \qu \text{for $t\in[0,T]$.} \label{lt}
\end{align}
\end{subequations}
By Implicit Function theorem and regularity of $(\rho,u)$ from Lemma \ref{lemma:local}, the derivatives $R_x$ and $R_t$ exist, and it follows from \eqref{nse1} and \eqref{Rexist} that for all $(x,t)\in\mL(T)$:
\begin{equation}\label{Rdiff}
	R_t(x,t) = u(R(x,t),t), \qu R_x(x,t) = \dfrac{R(x,t)^{1-n}}{\rho(R(x,t),t)}
\end{equation}
Therefore $R(x,t)$ satisfies the following identity:
\begin{align}\label{R}
&R(x,t)= R_0(x) + \int_{0}^{t} u(R(x,s),s)\dif s = \Big(1+n\int_{B(t)}^{x} \dfrac{1}{\rho(R(y,t),t)}\dif y\Big)^{\frac{1}{n}},\\
&\text{where $R_0(x)$ is implicitly defined by } \ x = \int_{1}^{R_0(x)} \rho_m^0(r) r^{n-1} \dif r \ \text{ for $0 \le x \le M_0$.}\nonumber
\end{align} 
The transformation from the Eulerian coordinate $(r,t)$ to the Lagrangian coordinate $(x,t)$ is executed by the equation
\begin{equation}\label{rR}
	r = R(x,t).
\end{equation}
Let $v \vcentcolon= 1/\rho$ be the specific volume. Using \eqref{Rdiff}, we obtain the system of equations for $(\hat{u}, \hat{v})(x,t):=(u,v)(R(x,t), t)$ from \eqref{nse} as
\begin{subequations}
\label{nsl-hat}
\begin{gather}
\hv_t = (r^{n-1} \hu)_x, \label{nsl-hat1} \\
\hu_t = \mu r^{n-1} \Big( \frac{(r^{n-1} \hu)_x}{\hv} \Big)_x
- r^{n-1} p(\hv)_x, \label{nsl-hat2}
\end{gather}
\end{subequations}
where $p(v)=K v^{-\gamma}$ and $r=R(x,t)$. The initial and the boundary conditions
for $(\hv,\hu)$ are derived from \eqref{approx2} and \eqref{approx3} as
\begin{subequations}\label{lc}
\begin{gather}
\hv(x,0) = \hv_0(x) := 1/\rho_m^0(R_0(x)), \qu \hu(x,0) = \hu_0(x) := u_m^0(R_0(x)), \label{lic-hat} \\
\hv_0(M_0)= \vt(m), \qu \hu_0(M_0,t) = \ut(m), \qu \hu_0(0)=u_b, \label{lcompa-hat}\\
\hu(B(t),t) = u_b, \qu \hu(M(t),t)= \ut(m), \qu \hv(M(t),t)=\vt(m)\vcentcolon=(\rt(m))^{-1}. \label{lbc-hat}
\end{gather}
\end{subequations}
Since the spatial variable $r = R(x,t)$ in the Eulerian coordinate depends on the spatial and time variables $(x,t)$ in the Lagrangian coordinate, the density $\tilde{\rho}$ of the stationary solution also depends on $(x,t)$. For simplicity we denote $\tilde{\rho}(x,t):= \tilde{\rho}(R(x,t))$. Consequently, the specific volume $\vt$
in the stationary solution is also a function of $(x,t)$. Namely, $\vt(x,t) := 1/\rt(R(x,t))$. In addition, $\tilde{\rho}_0(x):= \tilde{\rho}(R_0(x))$ and  $\tilde{v}_0(x):=1/\tilde{\rho}_0(x)$.

By \eqref{Rbdry}, we have $v(1,t)=v(R(B(t),t),t)=\hv(B(t),t)$ for each $t\in[0,T]$. Therefore the free boundary $B(t)$ can be rewritten as
\begin{equation*}
	B(t)=|u_b| \int_{0}^{t} \rho(1,s) \dif s = \int_{0}^{t} \dfrac{|u_b|}{\hv(B(s),s)}\dif s.
\end{equation*}
The condition \eqref{lbc-hat} along with chain rules imply that for each $t\in[0,T]$,
\begin{equation*}
	\hu_t(B(t),t) + B^{\prime}(t) \hu_x(B(t),t) =0, \qu \text{where } \ B^{\prime}(t)= \dfrac{|u_b|}{\hv(B(t),t)}.
\end{equation*} 
Thus 
the 
compatibility condition \eqref{compa} at $r=1$ is written in Lagrangian coordinates as
\begin{equation}\label{lcompa1-hat}
\hu_0\vert_{x=0}=u_b, \qu \Big\{ |u_b| \dfrac{(\hu_0)_x}{\hv_0} + \mu \Big( \dfrac{(\hu_0)_x}{\hv_0} \Big)_x - p(\hv_0)_x \Big\} \Big\vert_{x=0} = 0.
\end{equation}
In addition, we define 
\begin{equation*}
	\hat{\phi}(x,t)\vcentcolon= \hv(x,t)-\vt(x,t), \qu \hat{\psi}(x,t)\vcentcolon= \hu(x,t)-\ut(x,t).
\end{equation*}
Boundary condition \eqref{lbc} and \eqref{Rbdry} then imply that:
\begin{equation}\label{ldbc}
	\hat{\psi}(B(t),t)=0, \qu \hat{\psi}(M(t),t)=\hat{\phi}(M(t),t)=0, \qu \text{for all } \ t\in[0,T].
\end{equation}

\section{A-priori Estimates in Lagrangian Coordinates}\label{sec:apriori}\setcounter{equation}{0}
Until the end of Section \ref{sec:apriori},
we omit the hat ``\,$\hat{ \ }$\,'' to express 
functions in the Lagrangian coordinate. Thus, we have
\begin{subequations}
	\label{nsl}
	\begin{gather}
		v_t = (r^{n-1} u)_x, \label{nsl1} \\
		u_t = \mu r^{n-1} \Big( \frac{(r^{n-1} u)_x}{v} \Big)_x
		- r^{n-1} p(v)_x, \label{nsl2}
	\end{gather}
\end{subequations}
where $r=R(x,t)$. The initial and boundary data are given by
\begin{subequations}\label{libc}
\begin{gather}
	v(x,0) = v_0(x), \qu
	u(x,0) = u_0(x), \label{lic} \\
	u(B(t),t) = u_b, \qu v(M(t),t)=\vt(m), \qu u(M(t),t)=\ut(m),
	\label{lbc}
\end{gather}
\end{subequations}
and the initial data satisfies the compatibility condition:
\begin{subequations}\label{lcompa}
\begin{gather}
	u_0\vert_{x=0}=u_b, \qu \Big\{ |u_b| \dfrac{(u_0)_x}{v_0} + \mu \Big( \dfrac{(u_0)_x}{v_0} \Big)_x - p(v_0)_x \Big\} \Big\vert_{x=0} = 0,\label{lcompaIN}\\
	v_0(M_0)= \vt(m), \qu u_0(M_0,t) = \ut(m). \label{lcompaOUT}
\end{gather}
\end{subequations}
We consider the initial boundary value problem to the system of equations
\eqref{nsl} with data \eqref{lic} and \eqref{lbc}.
The coefficients in \eqref{nsl} 
is given by the relation \eqref{R} and \eqref{rR}.

We  distinguish the differential operators
in two different coordinate systems as follows: the spatial and temporal derivatives in Lagrangian coordinates are denoted respectively as $D_x$ and $D_t$, while the spatial and temporal derivatives in Eulerian coordinates are denoted respectively as $\partial_r$ and $\partial_t$. Using relations (\ref{Rdiff}), it can be verified that
\begin{equation*}
\begin{aligned}
D_x = \dfrac{v}{r^m} \partial_r, \quad D_t = \partial_t + u \partial_r, \quad \partial_r = \dfrac{r^m}{v} D_x, \quad \partial_t = D_t - \dfrac{r^m u}{v} D_x. 
\end{aligned}
\end{equation*}
In what follows, we   use
the notation $\partial_r f$ and $\partial_t f$ to indicate the Eulerian derivative of $f$, while we abbreviate as $f_x \equiv D_x f$ and $f_t \equiv D_t f$ for the Lagrangian derivatives. 

%The a-priori estimates of the difference functions $(\psi, \phi)$ for this problem is stated in the following theorem.

\subsection{Basic energy estimate}

In this subsection,
we show the basic energy estimate for the solution %
$(v, u)$ to the problem \eqref{nsl}--\eqref{lcompa} uniformly in $m$. 
To this end, we employ the energy
form
$\vE$ defined by
\begin{gather}
 \vE := \frac{1}{2} |u-\ut|^2 + G(v, \vt), \quad
\text{where} \qu G(v, \vt) := \int_{\vt^{-1}}^{v^{-1}} \dfrac{p(z^{-1})-p(\vt)}{z^2}\, \dif z.%=K \Big\{\dfrac{v^{1-\gamma}-\vt^{1-\gamma}}{\gamma-1}+\vt^{-\gamma}(v-\vt)\Big\}.
\label{G}
\end{gather}
For $\gamma\ge1$, and $v$, $\vt\ge 0$, $G(v,\vt)$ is a positive convex function such that $G(v,v)=0$ for $v> 0$.
Using \eqref{st}, %and $p(v)=Kv^{-\gamma}$, 
we have
\begin{equation}
	\partial_{v} G(v,\vt) = p(\vt)-p(v), \ \ \partial_{\vt} G(v,\vt) = \phi p^{\prime}(\vt)\vcentcolon= (v-\vt) p^{\prime}(\vt).
\end{equation}
$G(v,\vt)$ is also  rewritten as
\begin{equation}\label{g}
G(v, \vt) = \vt p(\vt) g \big( \frac{v}{\vt} \big), \quad \text{where } \ \  g(s) := s - 1 - \int_1^s \eta^{-\gamma} \, d \eta.
\end{equation}
The function $g(s)$ satisfies the lower bound:
\begin{equation}
s - 1 - \log s \le g(s) \quad \text{for} \ \ s > 0.
\label{psi}
\end{equation}
By (\ref{stProp2}), there exists 
a constant 
$\delta=\delta(\rho_+,\mu,\gamma,K,n)>0$ such that if $|u_b|\le \delta$, then the following point-wise upper and lower bounds hold for $\vt(x,t) \equiv \vt(R(x,t)) = 1/\rt(R(x,t))$: 
\begin{equation}\label{vtbds}
\dfrac{1}{\rho_+} \le \vt(x,t) \le \dfrac{2}{\rho_+}, \qu \text{for} \qu (x,t)\in\mL(T),
\end{equation}
\begin{proposition}\label{prop:phiG}
If $1\le\gamma\le 2$, then for $v$, $\vt\ge 0$,
\begin{equation*}
\dfrac{\gamma K}{2}|v-\vt|^2 \le \left\{ \begin{aligned}
&\vt^{1+\gamma} G(v,\vt) && \text{if $v\le \vt$,}\\
&\vt^{\gamma} v G(v,\vt) && \text{if $v>\vt$.}
\end{aligned} \right.
\end{equation*}
%If $\gamma=1$, then for $v$, $\vt\ge 0$,
%\begin{equation*}
%\dfrac{K}{2}|v-\vt|^2 \le \left\{ \begin{aligned}
%&\vt^{2} G(v,\vt) && \text{if $v\le \vt$,}\\
%&\vt v G(v,\vt) && \text{if $v>\vt$.}
%\end{aligned} \right.
%\end{equation*}
\end{proposition}
The proof for this proposition is rather technical and it is not the main point of this subsection. Hence it is placed in Appendix \ref{appen:G}.

%At first, we show the energy equality.
By the standard discussion using the mollifier with respect to time $t$,
we may regard the solution $(v,u) \in X_{\ul{d},\bar{d},D}^m(0,T)$ as if
it verifies
\begin{gather*}
(v - \vt,u) \in C^\infty([0,T]; H^1(0,m)),
\quad
u \in C^\infty((0,T]; H^2(0,m)).
\end{gather*}
%% basic est.
%%
\begin{lemma}\label{lemma:be}
Suppose $1\le \gamma\le 2$, and there exists $T>0$ such that $(v, u)$ solves \eqref{nsl}--\eqref{lcompa} in $\mL(T)$. Then there exists $\delta=\delta(\rho_+,\mu,\gamma,K,n)>0$ such that if $|u_b|\le \delta$ then
\begin{multline*}
\sup\limits_{t\in[0,T]}\int_{S(t)}\!\!\! \vE(x,t) \, \dif x
+ \mu\!\! \iint_{\mL(T)} \!\! \Big\{\dfrac{n-1}{2}\frac{v\psi^2}{r^2}
+ \frac{r^{2n-2}}{v} \psi_x^2 \Big\}  \, \dif x \dif t+|u_b|\!\!\int_{0}^{T} \dfrac{G(v,\vt)}{v} (B(t),t)\, \dif t\\
+C^{-1}\iint_{\mL(T)}\!\! \Big\{ |u_b|\dfrac{\psi^2}{r^n} + |u_b|^3\dfrac{G(v,\vt)}{r^{3n-2}} \Big\}  \, \dif x \dif t \le \int_0^{M_0}\!\!\!  \vE(x,0)\,  \dif x,
\end{multline*}
where $C=C(\rho_+,\mu,\gamma,K,n)>0$ is a constant independent of $m$, $T$.
\end{lemma}
\begin{proof}
Taking Lagrangian temporal derivative on $\mathcal{E}$, using the differential relations \eqref{Rdiff} and the equations for stationary solution \eqref{st}, we obtain that
\begin{multline*}
\vE_t + \mu \dfrac{(r^{n-1}\psi)_x^2}{v} + (\gamma-1) \rt(1) |u_b| \dfrac{\partial_r \rt }{r^{n-1}\rt^2} G(v,\vt) + \partial_r \ut |\psi|^2 \\
= \Big\{ \big( \mu \dfrac{(r^{n-1}\psi)_x}{v} + p(\vt)-p(v) \big) r^{n-1}\psi  \Big\}_x + \tilde{L} \psi \phi,
\end{multline*}
where $\tilde{L} \vcentcolon= \mu \partial_r ( r^{1-n} \partial_r (r^{n-1}\ut) )$. By (\ref{stProp1}), we have $\tilde{L} = \mu u_b \rt(1) \partial_r ( r^{1-n} \partial_r \vt )$. Integrating the above equality in $x\in S(t)\vcentcolon=[B(t),M(t)]$, we obtain with the help of 
the
boundary condition \eqref{lbc} that
\begin{multline}\label{beTemp1}
\int_{S(t)}\!\! \vE_t\, \dif x + \int_{S(t)}\!\! \Big\{ \mu \dfrac{(r^{n-1}\psi)_x^2}{v} + (\gamma-1) \rt(1) |u_b| \dfrac{\partial_r \rt }{r^{n-1}\rt^2} G(v,\vt) + \partial_r \ut \psi^2 \Big\}\, \dif x \\
= \int_{S(t)} \tilde{L} \psi \phi\, \dif x.
\end{multline}
The 
condition (\ref{lbc}) also implies that $\vE(M(t),t)=G(\tilde{v}(m),\tilde{v}(m))=0$. Thus applying Leibniz's Integral Rule and (\ref{B}), it follows that
\begin{align}\label{beTemp2}
\dfrac{\dif }{\dif t} \int_{S(t)} \vE (x,t) \dif x =& \int_{S(t)} \vE_t(x,t) \dif x + M^{\prime}(t) \vE(M(t),t) - B^{\prime}(t) \vE(B(t),t)\nonumber\\
=& - |u_b| \dfrac{G(v,\vt)}{v}(B(t),t) + \int_{S(t)} \vE_t(x,t) \dif x .
\end{align}
In addition, using 
the
boundary condition (\ref{lbc})  again, we obtain with integration by parts that
\begin{align}\label{beTemp3}
\int_{S(t)} \dfrac{(r^{n-1}\psi)_x^2}{v}\, \dif x %= \int_{S(t)} v^{-1} \big( r^{n-1} \psi_x + (n-1) \dfrac{v}{r} \psi  \big)^2 \dif x
=& \int_{S(t)} \big( \dfrac{r^{2(n-1)}}{v} \psi_x^2 + 2(n-1) r^{n-2} \psi_x \psi + (n-1)^2 \dfrac{v}{r^2} \psi^2  \big)\, \dif x\nonumber \\
=& \int_{S(t)} \big( \dfrac{r^{2(n-1)}}{v} \psi_x^2 + (n-1) \dfrac{v}{r^2} \psi^2  \big)\, \dif x. 
\end{align}
Substituting (\ref{beTemp2})--(\ref{beTemp3}) in (\ref{beTemp1})
yields that
\begin{align}\label{beTemp4}
& \dfrac{\dif}{\dif t}\int_{S(t)}\!\! \vE(x,t)\, \dif x + |u_b| \dfrac{G(v,\vt)}{v}(B(t),t) + \mu \int_{S(t)} \!\!\Big\{ \dfrac{r^{2n-2}}{v} \psi_x^2 + (n-1) \dfrac{v}{r^2} \psi^2 \Big\}(x,t)\, \dif x \nonumber\\
&+\int_{S(t)}\!\!\Big\{ (\gamma-1)\rt(1)|u_b|\dfrac{\partial_r \rt}{r^{n-1}\rt^2 }G(v,\vt) + \partial_r \ut \psi^2 \Big\}(x,t)\, \dif x = \int_{S(t)} \tilde{L} \phi \psi (x,t) \, \dif x. 
\end{align}
From Lemma \ref{lemma:st}, the following estimates holds for the stationary solution:
\begin{align*}
&(\gamma-1)\rt(1)|u_b|\dfrac{\partial_r \rt}{r^{n-1}\rt^2} \ge C_1^{-1} |u_b|^3 r^{-3n+2}, \qquad \partial_r \ut \ge C_2^{-1} |u_b| r^{-n},\\
&|\tilde{L}| %= \mu \rt(1) |  u_b |\cdot  \big| \partial_r \big( \dfrac{\partial_r \rt}{r^{n-1}\rt^2} \big) \big|
= \mu \rt(1) |  u_b |\cdot \big| \dfrac{\partial_r^2 \rt}{r^{n-1}\rt^2} -(n-1) \dfrac{\partial_r \rt }{r^{n}\rt^2} - 2 \dfrac{|\partial_r \rt|^2}{r^{n-1}\rt^3} \big| \le C_3 |u_b|^3 r^{-3n+1},
\end{align*}
where for each $i=1, 2, 3$,  $C_{i}=C_{i}(\rho_+,\mu,\gamma,K,n)>0$ are some constants. 
Substituting these inequalities in (\ref{beTemp4}), we have that
\begin{align}\label{beTemp5}
&\dfrac{\dif}{\dif t}\int_{S(t)}\!\! \vE(x,t)\, \dif x + |u_b| \dfrac{G(v,\vt)}{v}(B(t),t) + \mu \int_{S(t)} \!\!\Big\{ \dfrac{r^{2n-2}}{v} \psi_x^2 + (n-1) \dfrac{v}{r^2} \psi^2 \Big\}(x,t)\, \dif x \nonumber\\
&+\int_{S(t)}\!\!\Big\{ \dfrac{|u_b|^3}{C_1} \dfrac{G(v,\vt)}{r^{3n-2}} + \dfrac{|u_b|}{C_2} \dfrac{\psi^2}{r^n} \Big\} (x,t) \, \dif x \le C_3 |u_b|^3 \int_{S(t)} \dfrac{|\phi |\cdot|\psi|}{r^{3n-1}}(x,t)\, \dif x. 
\end{align}
Owing to
Proposition \ref{prop:phiG}, we estimate the second term on the right hand side of (\ref{beTemp5}) in two cases: $\{ v\le \vt \}$ and $\{ v > \vt \}$. First we consider $\{v\le \vt\}$. Then by the Cauchy-Schwartz inequality and (\ref{vtbds}), it holds that
\begin{align}\label{beTemp6}
&C_3 |u_b|^3 \int_{S(t)\cap \{ v \le \vt \}} \dfrac{|\phi| \cdot |\psi| }{r^{3n-1}}\, \dif x \le C_3|u_b|^3 \sqrt{\dfrac{2}{\gamma K}} \int_{S(t)\cap \{ v \le \vt \}} \dfrac{|\psi|}{r^{3n-1}}  \sqrt{\vt^{1+\gamma}G(v,\vt)}\, \dif x \nonumber\\
&\le |u_b|^5\dfrac{2^{1+\gamma}C_3^2 C_2}{\gamma K \rho_+^{1+\gamma}} \int_{S(t)} \dfrac{G(v,\vt)}{r^{5n-2}}\, \dif x + \dfrac{|u_b|}{2 C_2} \int_{S(t)} \dfrac{\psi^2}{r^n}\, \dif x,\nonumber\\
&\le \dfrac{|u_b|^3}{4 C_1}  \int_{S(t)} \dfrac{G(v,\vt)}{r^{3n-2}}\, \dif x + \dfrac{|u_b|}{2C_2} \int_{S(t)} \dfrac{\psi^2}{r^n}\, \dif x,
\end{align}
where in the last line we have chosen $|u_b|\le \delta_1$ with
\begin{equation*}\label{delta1}
\delta_1 \vcentcolon= \sqrt{\dfrac{\gamma K \rho_+^{1+\gamma}}{ 2^{3+\gamma} C_3^2 C_1 C_2}}.
\end{equation*}
Moreover, we  have the following estimate for $\{ v > \vt \}$ using Proposition \ref{prop:phiG}:
\begin{align}\label{beTemp7}
&C_3 |u_b|^3 \int_{S(t)\cap \{v > \vt\}} \dfrac{|\phi|\cdot |\psi|}{r^{3n-1}}\, \dif x  \le C_3 |u_b|^3 \sqrt{\dfrac{2}{\gamma K}} \int_{S(t)\cap \{ v > \vt \}} \dfrac{|\psi|}{r^{3n-1}} \sqrt{\vt^{\gamma} v G(v,\vt)}\, \dif x \nonumber\\
&\le \dfrac{(n-1)\mu}{2} \int_{S(t)} \dfrac{v \psi^2}{r^2 }\, \dif x + |u_b|^6 \dfrac{2^{\gamma} C_3^2}{(n-1)\mu \gamma K \rho_+^{\gamma}}\int_{S(t)} \dfrac{G(v,\vt)}{r^{6n-4}}\, \dif x,\nonumber\\
&\le \dfrac{(n-1)\mu}{2} \int_{S(t)} \dfrac{v \psi^2}{r^2 }\, \dif x +  \dfrac{|u_b|^3}{4C_1} \int_{S(t)} \dfrac{G(v,\vt)}{r^{6n-4}}\, \dif x,
\end{align}
where we have chosen $|u_b|\le \delta_2$ with 
\begin{equation*}
    \delta_2 \vcentcolon= \sqrt[\leftroot{-1}\uproot{2}\scriptstyle 3]{\frac{(n-1)\mu \gamma K \rho_+^{\gamma} }{2^{2+\gamma} C_3^2 C_1 }}.
\end{equation*}
Substituting (\ref{beTemp6})--(\ref{beTemp7}) in (\ref{beTemp5}), it holds that for $|u_b|\le \delta \vcentcolon= \min\{\delta_1,\delta_2\}$, 
\begin{align}\label{beConcl}
&\dfrac{\dif }{\dif t}\int_{S(t)}\!\! \vE(x,t)\, \dif x + |u_b| \dfrac{G(v,\vt)}{v}(B(t),t) + \mu \int_{S(t)} \!\!\Big\{ \dfrac{r^{2n-2}}{v} \psi_x^2 + \dfrac{n-1}{2} \dfrac{v \psi^2}{r^2} \Big\}(x,t)\, \dif x \nonumber\\
&+\int_{S(t)}\!\!\Big\{ \dfrac{|u_b|^3}{2C_1} \dfrac{G(v,\vt)}{r^{3n-2}} + \dfrac{|u_b|}{2C_2} \dfrac{\psi^2}{r^n} \Big\}(x,t)\, \dif x \le 0.
\end{align}
Integrating the above inequality in $[0,t]$, then taking the supremum over $t\in[0,T]$, we obtain the desired estimate.
\end{proof}

\bigskip

%% corollary
%%
The following lemma is proved by the Sobolev inequality.
It is utilized in deriving the pointwise estimate in
Subsection \ref{subsec:point}.
\begin{corollary}\label{corol:infpsi}
Suppose assumptions of Lemma \ref{lemma:be} holds, then
\begin{equation*}
 \int_0^T \sup\limits_{ x\in S(t)}|r^{n-2} \psi^2|(x,t) \, \dif t
\le \int_0^{M_0} \vE(x,0) \, \dif x.
\end{equation*}
\end{corollary}
\begin{proof}
Fix a time $t\in[0,T]$. Since \eqref{lbc}, for all $x\in [B(t),M(t)]$,
\begin{equation*}
(r^{n-2} \psi^2)(x,t)
=
\int_{B(t)}^x (r^{n-2} \psi^2)_x(x,t) \, \dif x
=
\int_{B(t)}^x \{ (n-2) r^{n-3} r_x \psi^2 + 2 r^{n-2} \psi \psi_x \}(x,t) \, \dif x.
%\label{b2}
\end{equation*}
Substituting $r_x = v/r^{n-1}$ in the above due to \eqref{Rdiff}, and
taking the absolute value of the resulting equality
and applying the Schwarz inequality, we have
\begin{equation*}
\sup\limits_{x\in S(t)}| r^{n-2} \psi^2| (x,t)
 \le \int_{B(t)}^{M(t)} \Big\{  (n-1) \frac{v}{r^2} \psi^2 + \frac{r^{2n-2}}{v} \psi_x^2 \Big\}(x,t) \, \dif x.
\label{b3}
\end{equation*}
Integrating the above inequality over $(0,T)$, then using Lemma \ref{lemma:be}, we obtain the desired estimate.
\end{proof}

\begin{corollary}\label{corol:Bbound}
Suppose assumptions of Lemma \ref{lemma:be} holds. Then
\begin{equation}
B(T) \le C_0 ( 1+ |u_b|T ),
\end{equation}
where $C_0=C_0(\rho_0,u_0,\rho_+,\mu,\gamma,K,n)>0$ is a constant independent of $m$ and $T$.
\end{corollary}
\begin{proof}
To avoid confusion, the notations used in Section \ref{sec:RL} will be reinstated for the proof of Corollary \ref{corol:Bbound}. More specifically, 
$\rho(r,t)$ denotes the density in Eulerian coordinate, whereas $\hat{\rho}(x,t)\equiv \rho(R(x,t),t)$  denotes 
the density in Lagrangian coordinate. 

For $\rho$, $\rt\ge 0$, we set $H(\rho,\rt)\vcentcolon= \rho G(\rho^{-1},\rt^{-1})$. Then it follows that
\begin{equation}\label{Hh}
H(\rho,\rt) = \left\{ \begin{aligned} &\dfrac{K\rt^{\gamma}}{\gamma-1}  h\big(\dfrac{\rho}{\rt}\big), && \text{where } \ h(s) = s^{\gamma} - 1 -\gamma (s-1) && \text{if} \qu \gamma>1,\\
& K \rt h\big(\dfrac{\rho}{\rt}\big), && \text{where } \ h(s)= s\log s - s + 1 && \text{if} \qu \gamma=1. 
\end{aligned} \right.
\end{equation}
For both cases $\gamma>1$ or $\gamma=1$, it is verified that the map $s\mapsto h(s)$ is non-negative, smooth, and convex in $[0,\infty)$, and it achieves a unique absolute minimum at $s=1$ with $h(1)=0$. By the boundary condition (\ref{lbc}), we have $\hat{v}(B(t),t)=1/\rho(1,t)$. Thus 
Lemma \ref{lemma:be} gives the estimate:
\begin{equation}\label{bebdry}
|u_b|\int_{0}^T H(\rho(1,t),\rt(1)) \dif t = |u_b|\int_{0}^{T} \dfrac{G(\hat{v},\vt)}{\hat{v}}(B(t),t) \dif t \le E_0.
\end{equation}
Since $h(s)$ is convex, it follows by Jensen's inequality, (\ref{vtbds}), and (\ref{Hh}) that
\begin{equation*}
h\Big( \dfrac{1}{T} \int_{0}^{T} \dfrac{\rho(1,t)}{\rt(1)} \dif t \Big) \le \dfrac{1}{T} \int_0^T h\big(\dfrac{\rho(1,t)}{\rt(1)}\big) \dif t \le  \dfrac{\tilde{C}_0}{|u_b|T} \vcentcolon=\left\{\begin{aligned}
&\dfrac{2^{\gamma}(\gamma-1)E_0}{K\rho_+^{\gamma}}\cdot \dfrac{1}{|u_b|T} && \text{if } \ \gamma>1,\\
&\dfrac{2 E_0}{K\rho_+} \cdot \dfrac{1}{|u_b|T} && \text{if } \ \gamma=1.
\end{aligned}   \right.
\end{equation*}
Let $q\vcentcolon [0,\infty)\to [1,\infty)$ be the right branch inverse of $h(\cdot)$. Then above inequality implies that
\begin{equation}\label{Bineq}
\int_{0}^{T} \rho(1,t) \dif t \le \dfrac{\rt(1)}{|u_b|} |u_b|T q\big(\dfrac{\tilde{C}_0}{|u_b|T}\big) \le  \dfrac{\rho_+}{|u_b|} |u_b|T q\big(\dfrac{\tilde{C}_0}{|u_b|T}\big).
\end{equation}
Since $q$ is the inverse map of $h$, it is shown by Implicit Function theorem that if $\gamma>1$ then $q^{\prime}(z)=\gamma^{-1}(q(z)^{\gamma-1}-1)^{-1}$, and if $\gamma=1$ then $q^{\prime}(z)=(\log q(z))^{-1}$ for $z\in(0,\infty)$. Moreover we also have $q(z)\to\infty$ as $z\to\infty$. Thus by L'Hospital's rule, $y q(\tilde{C}_0/y) \to 0$ as $y\to 0^+$. It then follows that there exists $\alpha_0=\alpha_0(\tilde{C}_0)>0$ independent of $|u_b|>0$, $T>0$, and $m\in\mathbb{N}$ such that $y q(\tilde{C}_0/y)\le 1$ for all $y\le \alpha_0$. Furthermore, $q(\tilde{C}_0/y )\le q(\tilde{C}_0/\alpha_0)$ for all $ y \ge \alpha_0$ since $z\mapsto q(z)$ is strictly monotone increasing. Combining the two cases, $y<\alpha_0$ and $y\ge \alpha_0$, we obtain that
\begin{equation*}
B(t)= |u_b|\int_{0}^T \rho(1,t)\dif t \le \rho_+ |u_b|T q\big( \dfrac{\tilde{C}_0}{|u_b| T} \big) \le C_0 (1+|u_b|T), 
\end{equation*}
where $C_0 \vcentcolon= \rho_+ \cdot \max \{1,q(\tilde{C}_0/\alpha_0)\}$. This is the desired estimate for $B(t)$.
\end{proof}

We recall from (\ref{MODE}) that $M_0=M_0(m)\to \infty$ as $m\to\infty$. Combining with Corollary \ref{corol:Bbound}, it follows that for 
sufficiently 
large integer $m\in\mathbb{N}$, 
there exists  $T>1$, an existence time, which is independent of $m$ such that
\begin{equation}\label{BTM0}
    B(T) \le M_0 - 2 
\end{equation}

Throughout the rest of the paper, we denote $C>0$ and $C_0>0$ to be positive constants with the following dependence:
\begin{equation}\label{constants}
C=C(\rho_+,\mu,\gamma,K,n), \qu \text{and} \qu C_0=C_0(\rho_0,u_0,\rho_+,\mu,\gamma,K,n).
\end{equation}
\bigskip

\subsection{Point-wise estimate of the specific volume}\label{subsec:point}
This subsection is devoted to the point-wise upper and lower bounds of specific volume $v(x,t)$ uniformly
in $m\in\mathbb{N}$ and $T>0$.
\begin{lemma}\label{lemma:v}
Suppose there is a time $T>0$ such that $(v,u)$ is the solution to the problem \eqref{nsl}--\eqref{lcompa} in $\mL(T)$. Then there exists $\delta_0=\delta_0(\rho_0,u_0,\rho_+,\mu,\gamma,K,n)>0$ such that if $|u_b|\le \delta_0$, then the specific volume $v$ satisfies
\begin{equation}\label{vulb}
\underline{v} \le v(x,t) \le \ol{v} \ \  \text{ for } \ (x,t)\in\mL(T), 
\end{equation}
where $0<\underline{v}<\ol{v}<\infty$ are positive constants with the dependence $\underline{v}=\underline{v}(C_0)$ and $\overline{v}=\overline{v}(C_0)$.
\end{lemma}
The proof of Lemma \ref{lemma:v} is divided into
the several lemmas and completed at the end of this subsection.
The next lemma follows from Lemma \ref{lemma:be}.
%% lemma : positivity of integration of v
%%
\begin{lemma}\label{lemma:vInt}
Suppose that the same assumptions as in Lemma \ref{lemma:v} hold. Then, there exist a constant $C_0>0$ independent of $m$ and $T$ such that for each $t\in[0,T]$,
\begin{equation}
C_0^{-1} \le \int_{I(t)} v(x,t) \, \dif x \le C_0, \ \  \text{ for any set $I(t)\subseteq S(t)$ satisfying $|I(t)|=1$.} \label{c1}
\end{equation}
\end{lemma}
\begin{proof}
Due to Lemmas \ref{lemma:st} and \ref{lemma:be}, there exists a constant $C_0>0$ such that
\begin{equation}\label{c2}
\int_{I(t)} g \big( \frac{v}{\vt} \big)(x,t) \, \dif x \le C_0,
\end{equation}
where $C_0>0$ is specified in (\ref{constants}). It is easy to see from (\ref{g}) that $g(s)$ satisfies $g^{\prime\prime}(s)>0$ for $s>0$, $g(1) = 0$,
	and $g(s) \to + \infty$ as $s \to + \infty$ or $s \to +0$. Thus, the equation $g(s)=C_0$ has the distinct two roots $\alpha_0$ and $\beta_0$
	satisfying
	$0 < \alpha_{0} < 1 < \beta_{0}$.
	Applying Jensen's inequality to the convex function $g(s)$
	and using (\ref{c2}), we have
	\begin{equation}\label{c4}
		g \Big( \int_{I(t)} \frac{v}{\vt} (x,t) \, \dif x \Big)
		\le \int_{I(t)} g \big( \frac{v}{\vt} \big)(x,t) \, \dif x
		\le C_0.
	\end{equation}
	Thus, we see that
	\begin{equation}\label{c41}
		\alpha_{0}   \le \int_{I(t)} \frac{v}{\vt} (x,t) \, \dif x
		\le \beta_{0},
	\end{equation}
 where $\alpha_{0}$ and $\beta_{0}$ are positive constants
 independent of $m$ and $T$.
	Combining
	the estimate
	(\ref{c41}) with (\ref{vtbds}) gives the desired inequality (\ref{c1}).
\end{proof}

\bigskip

\bigskip

%% representation formula
%%
Due to the moving boundary curves $B(t)$ and $M(t)$, derivation for the point-wise estimate of the specific volume $v(x,t)$ is conducted separately in 3 sub-regions of the Lagrangian space-time domain $\mL(T)$, which are:
\begin{subequations}\label{2regions}
\begin{align}
\mL_1(T) \vcentcolon=& \{ (x,t) \;|\; t\in[0,T], \ B(t)\le x \le M_0-2 \}, \label{L1}\\
\mL_2(T) \vcentcolon=& \{ (x,t) \;|\; t\in[0,T], \ M_0 \le x \le M(t) \},\label{L2}\\
\mL_3(T) \vcentcolon=& \{(x,t) \;|\; t\in[0,T], \ M_0-2 < x < M_0 \}.\label{L3}
\end{align} 
\end{subequations}
Note that $\mL_1(T)$ is well defined due to (\ref{BTM0}).

First, we consider the case that $(x,t)\in \mL_1(T)$. For each $t\in[0,T]$ and $k\in\mathbb{N}$, we employ the ``cut-off'' function $\zeta_{k,\,t}\vcentcolon [0,\infty) \to [0,1]$, given by
\begin{equation}\label{zeta}
	\zeta_{k,\,t}(y) :=
	\left\{
	\begin{array}{ll}
		1,     & \text{for} \ y\in [B(t) + k-1, B(t) + k], \\
		 1 - y + B(t) + k, & \text{for} \ y\in [B(t) + k, B(t)+k+1], \\
		0,     & \text{for} \ y\in [B(t) + k+1,\infty).
	\end{array}
	\right.
\end{equation}
In addition, for each $(k,t)\in\mathbb{N}\times[0,T]$, we define the closed interval
\begin{equation}\label{Jtk}
\mI_{t}^{k}\vcentcolon= [B(t)+k,B(t)+k+1].
\end{equation}
It follows that $\zeta_{k,t}$ satisfies the properties
\begin{equation}\label{zetaProp}
0\le \zeta_{k,\,t}(y) \le 1, \qu -1\le\zeta_{k,\,t}^{\prime}(y)\le 0, \qu \supp(\zeta_{k,\,t}^{\prime}) = \mI_{t}^k, \ \ \text{ for } \ y\in [0,\infty).
\end{equation}
We also define the function
\begin{equation}\label{alpha}
f(r)\vcentcolon= \rt(1) |u_b| \int_{1}^r \dfrac{\vt^{\prime}(s)}{s^{2(n-1)}} \dif s.
\end{equation}
\begin{lemma}\label{lemma:form1}
	Fix $(k,t)\in\mathbb{N}\times[0,T]$, and suppose assumptions in Lemma \ref{lemma:v} hold. Then
	\begin{equation}
		v(x,\tau)^{\gamma} =
		\frac{v_0(x)^\gamma + \frac{K \gamma}{\mu} \int_0^\tau A_{k,\,t}^1(x,s) D_{k,\,t}^1(x,s) \dif s}{A_{k,\,t}^1(x,\tau) D_{k,\,t}^1(x,\tau)},
		\label{form1}
	\end{equation}
	for $x\in \mI_t^{k-1}\cap[B(t),M_0-2]$ and $\tau\in[0,t]$, where $A_{k,\,t}^1(x,\tau)\vcentcolon=V^1_{k,\,t}(\tau) W^1_{k,\,t}(x,\tau)$ and
	\begin{alignat}{2}
		%A_{k,\,t}^1(x,\tau) &  := \exp \Big( && \frac{K \gamma}{\mu} \int_0^{\tau} \!\! \int_{\mI_{t}^k} v^{-\gamma} \, \dif y \dif s - \frac{(n-1)\gamma}{\mu} \int_0^{\tau} \!\! \int_x^{M_0} \zeta_{k,\,t} \frac{|\ut|^2}{r^{n}} \, \dif y \dif s \Big),  \\
		V^1_{k,\,t}(\tau) &\vcentcolon = \exp\Big(&& \dfrac{K\gamma}{\mu} \int_{0}^{\tau} \int_{\mI_{t}^k} v^{-\gamma} \dif y \dif \tau \Big), \label{V1} \\
		 W^1_{k,\,t}(x,\tau) &\vcentcolon= \exp\Big(&&-  \dfrac{(n-1)\gamma}{\mu} \int_{0}^{\tau} \int_{x}^{M_0} \zeta_{k,\,t} \dfrac{|\ut|^2}{r^{n}} \dif y \dif s \Big), \label{W1}\\
		D_{k,\,t}^1(x,\tau) & :=
		\exp \Big( && \frac{\gamma}{\mu} \int_x^{M_0} \zeta_{k,t}
		\big\{ \frac{\psi}{r^{n-1}} - f(r) \big\}(y,s) \, \dif y\Big\vert_{s=0}^{s=\tau} \nonumber \\
		& \ &&  + \frac{(n-1)\gamma}{\mu} \int_0^{\tau} \!\! \int_x^{M_0} \zeta_{k,t} \frac{\psi^2}{r^n} \, \dif y \dif s 
		- \gamma \int_{\mI_t^k} \log \frac{v(y,\tau)}{v_0(y)} \, \dif y
		\Big).\label{D1}
	\end{alignat}
\end{lemma}
\begin{proof}
    We consider the case $\mI_{t}^{k-1}\cap [B(t),M_0-2] \neq \emptyset$ as otherwise the statement is trivially true. Then, it follows that $\zeta_{k,\,t}(M_0)=0$ due to the definitions (\ref{zeta}) and (\ref{Jtk}). Multiplying (\ref{nsl2}) by $r^{1-n}$, rewriting the resulting equation with $u=\psi+\ut$, then using equations (\ref{st}), (\ref{stu}), and the relation (\ref{Rdiff}), we gets
	\begin{equation}\label{d4}
	    \Big( \dfrac{\psi}{r^{n-1}} - f(r) \Big)_s + (n-1) \dfrac{|\psi|^2}{r^n} - (n-1) \dfrac{|\ut|^2}{r^{n}} = \big( \mu (\log v)_s - K v^{-\gamma} \big)_y
	\end{equation}
	for $(y,s)\in\mL(T)$. Since $\rho>0$, the function
 $B(t)$ is strictly monotone increasing by (\ref{B}). Hence we see that for each fixed $t\in[0,T]$, if $B(t)\le x$ and $\tau\in[0,t]$ then
	\begin{equation*}
	\{ (y,s) \;|\; s\in[0,\tau], \text{ and }  x\le y\le M_0  \} \subset \mL(T).
	\end{equation*}
	Thus for $(x,\tau) \in \mI_{t}^{k-1}\times[0,t]$, we
 perform the following operations: multiplying (\ref{d4}) by $\zeta_{k,\,t}$ and
	integrating the resultant equality over $y\in [x, M_0]\times[0,\tau]$, then using the fact that $\zeta_{k,\,t}(M_0)=0$, we have
	\begin{align*}
		&\int_x^{M_0}\!\! \zeta_{k,t} 
		\big\{ \frac{\psi}{r^{n-1}}
		- f(r)
		\big\}(y,s) \, \dif y \Big\vert_{s=0}^{s=\tau}
		+ \int_0^{\tau} \!\! \int_x^{M_0}\!\! \zeta_{k,\,t} \dfrac{n-1}{r^n} ( |\psi|^2 - |\ut|^2 ) \, \dif y \dif s \nonumber\\
		&=
		K \int_0^{\tau} v(x,s)^{-\gamma} \, \dif s - \mu \log \frac{v(x,\tau)}{v_0(x)} + \mu \int_{\mI_t^k} \log \dfrac{v(y,\tau)}{v_0(y)}\dif y
		- K \int_0^{\tau} \!\! \int_{\mI_t^k}  v^{-\gamma} \, \dif y \dif s.
	\end{align*}
	Multiplying the above equality by $\gamma / \mu$
	and taking exponential of the resultant equality, we obtain that
	\begin{subequations}
		\begin{align}
			v_0(x)^{-\gamma} A_{k,\,t}^1(x,\tau) \, D_{k,\,t}^1(x,\tau)
			& =
			v(x,\tau)^{-\gamma}
			\exp \Big(
			\frac{K \gamma}{\mu} \int_0^{\tau} v(x,s)^{-\gamma} \, \dif s \Big)
			\label{d6} \\
			& =
			\frac{\mu}{K \gamma} \, \frac{\dif}{\dif \tau}
			\exp \Big(
			\frac{K \gamma}{\mu} \int_0^{\tau} v(x,s)^{-\gamma} \, \dif s \Big).
			\label{d7}
		\end{align}
	\end{subequations}
	Multiplying (\ref{d7}) by $K\gamma / \mu$, and integrating resultant equality
	over $[0,\tau]$, we have
	\begin{equation}
		\exp \Big(
		\frac{K \gamma}{\mu} \int_0^{\tau} v(x,s)^{-\gamma} \, \dif s
		\Big)
		=
		1 + v_0(x)^{-\gamma} \frac{K\gamma}{\mu}
		\int_0^{\tau} A_{k,\,t}^1(x,s) D_{k,\,t}^1(x, s) \, \dif s.
		\label{d8}
	\end{equation}
	Substituting (\ref{d8}) in (\ref{d6}) yields
	the desired formula (\ref{form1}).
\end{proof}

Next, we consider the case $(x,t)\in\mL_2(T)$. Set $t_x = (M_0-x)/(\rt(1)u_b)>0$, then $x=M(t_x)$ by (\ref{M}), and $v(x,t_x)=v(M(t_x),t_x)=\vt(m)$ by the boundary condition (\ref{lbc}). For $(k,t)\in\mathbb{N}\times[0,T]$, we define the cut-off function $\xi_{k,\,t}\vcentcolon [0,\infty)\to [0,1]$ as   
\begin{equation}\label{xi}
	\xi_{k,\,t}(y) \vcentcolon=
	\left\{
	\begin{array}{ll}
		0,     & \text{for} \ y\in[0,M(t)-k-1], \\
		y - M(t) +k +1, & \text{for} \ y\in [M(t) -k-1, M(t)-k], \\
		1,     & \text{for} \ y\in [M(t) - k, \infty).
	\end{array}
	\right.
\end{equation}
We also define the closed interval
\begin{equation}\label{mJtk}
\mJ_t^k\vcentcolon= [M(t)-k-1,M(t)-k], \qu \text{ for } \ (k,t)\in\mathbb{N}\times[0,T]
\end{equation}
Using these, we obtain the representation formula for the case $(x,t)\in\mL_2(T)$.
\begin{lemma}\label{lemma:form2}
Fix $(k,t)\in\mathbb{N}\times[0,T]$, and suppose assumptions in Lemma \ref{lemma:v} hold. Then
	\begin{equation*}
		v(x,\tau)^{\gamma} =
		\frac{\vt(m)^\gamma + \frac{K \gamma}{\mu} \int_{t_x}^\tau A^2_{k,\,t}(x,s) D^2_{k,\,t}(x,s) \dif s}{A^2_{k,\,t}(x,\tau) D^2_{k,\,t}(x,\tau)},
	\end{equation*}
	for all $x\in \mJ_t^{k-1}\cap[M_0,M(t)]$ and $\tau\in[t_x,t]$, where $A^2_{k,\,t}(x,\tau)\vcentcolon=V_{k,\,t}^2(\tau) W^2_{k,\,t}(x,\tau)$ and
	\begin{alignat*}{2}
		%A^2_{k,\,t}(x,\tau) &  := \exp \Big( && \frac{K \gamma}{\mu} \int_{t_x}^{\tau} \!\! \int_{\mJ_{t}^k} v^{-\gamma} \, \dif y \dif s + \frac{(n-1)\gamma}{\mu} \int_{t_x}^{\tau} \!\! \int_{B(s)}^{x} \xi_{k,\,t} \frac{|\ut|^2}{r^{n}} \, \dif y \dif s \Big),\\
		 V_{k,\,t}^2(\tau) &\vcentcolon = \exp\Big(&& \dfrac{K\gamma}{\mu} \int_{t_x}^{\tau} \int_{\mJ_{t}^k} v^{-\gamma} \dif y \dif \tau \Big),\\
		 W^2_{k,\,t}(x,\tau) &\vcentcolon= \exp\Big(&&  \dfrac{(n-1)\gamma}{\mu} \int_{t_x}^{\tau} \int_{B(s)}^{x} \xi_{k,\,t} \dfrac{|\ut|^2}{r^{n}} \dif y \dif s \Big),\\
		D^2_{k,\,t}(x,\tau) &\vcentcolon =
		\exp \Big( &&\frac{\gamma}{\mu} \int_{B(s)}^{x} \xi_{k,t}
		\big\{ \frac{\psi}{r^{n-1}} - f(r) \big\}(y,s) \, \dif y\Big\vert_{s=\tau}^{s=t_x} \\
		& \ &&  - \frac{(n-1)\gamma}{\mu} \int_{t_x}^{\tau} \!\! \int_{B(s)}^{x} \xi_{k,t} \frac{\psi^2}{r^n} \, \dif y \dif s 
		- \gamma \int_{\mJ_t^k} \log \frac{v(y,\tau)}{v(y,t_x)} \, \dif y
		\Big).
	\end{alignat*}
\end{lemma}
The proof of Lemma \ref{lemma:form2} follows almost identical procedure to that of Lemma \ref{lemma:form1}, except 
that 
we have used the cut-off function $\xi_{k,\,t}$, and the domain of integration is $(y,s)\in (-\infty,x]\times[t_x,\tau]$. 

\bigskip

%Using Lemmas \ref{lemma:form1}--\ref{lemma:form2}, we obtain the point-wise upper and lower bounds on $v(x,t)$ in the domain $\mL_1(T)\cup \mL_2(T)$, which are uniform over the parameters $m$, $T$. The derivations of this for $\mL_1(T)$ and $\mL_2(T)$ are almost identical. Thus, in the ensuing lemmas, we will only present the proof for the case $(x,t)\in\mL_1(T)$.

%% estimate of D(x,t)
%%

\begin{lemma}\label{lemma:D}
Fix $(k,t)\in\mathbb{N}\times[0,T]$. Suppose assumptions in Lemma \ref{lemma:v} hold. Then there exists a constant $C_0>0$ independent of $m$, $T$, $t$ and $k$, such that
\begin{subequations}
\begin{align}
&C_0^{-1} \le D^1_{k,\,t}(x,\tau) \le C_0 \ 
&&\text{ for all } x\in [B(t),M_0-2]\cap\mI_{t}^{k-1} \ \text{ and } \ \tau\in \times[0,t],\label{f0}\\
&C_0^{-1} \le D^2_{k,\,t}(x,\tau) \le C_0 \ 
&&\text{ for all }  x\in [M_0,M(t)]\cap \mJ_t^{k-1} \ \text{ and } \ \tau \in [t_x,t].	\label{f00}
\end{align}
\end{subequations}
\end{lemma}
\begin{proof}
We will only present the proof for (\ref{f0}) since (\ref{f00}) can be obtained by the same argument. %The desired estimate (\ref{f0}) follows from the estimation of each term in (\ref{D1}). 
Applying the Cauchy-Schwarz inequality and using Lemma \ref{lemma:be}, the first term of the right-hand side of (\ref{D1}) is estimated as
\begin{equation}
\Big| \frac{\gamma}{\mu} \int_x^{M_0} \zeta_{k,\,t}\big\{ \frac{\psi}{r^{n-1}} - f(r) \big\}(y,s) \, \dif y \Big|
\le C \int_{S(s)} \psi^2 (y,s) \, \dif y + C \int_x^{M_0} \zeta_{k,\,t} \, \dif y \le C_0, \label{f1}
\end{equation}
for $s=0$ or $\tau$. Using Lemma \ref{corol:infpsi}, the second term is estimated as
\begin{equation}\label{f2}
\Big| \frac{(n-1)\gamma}{\mu} \int_0^{\tau} \!\! \int_x^{M_0} \zeta_{k,\,t} \frac{|\psi|^2}{r^n} \, \dif y \dif s \Big| \le C \int_0^\tau \!\! \sup\limits_{y\in S(s)} \psi^2(y,s) \, \dif s \le C E_0.
\end{equation}
Applying Jensen's inequality on (\ref{c1}), we have the upper bound for the last term as
\begin{equation}\label{f3}
\int_{\mI_{t}^k} \log \frac{v(y,\tau)}{v_0(y)} \, \dif y \le \log \Big( \int_{\mI_{t}^k} \frac{v(y,\tau)}{v_0(y)} \, \dif y \Big) \le C_0.
\end{equation}
The lower bound  follows from Lemma \ref{lemma:be} and (\ref{c1}) as
\begin{align}\label{f4}
\int_{\mI_{t}^k} \log \frac{v(y,\tau)}{v_0(y)} \, \dif y & = \int_{\mI_{t}^k} \big\{ \log \frac{v(y,\tau)}{\vt(y,\tau)} + \log \frac{\vt(y,\tau)}{v_0(y)} \big\} \, \dif y \nonumber \\
& \ge - \int_{\mI_t^{k}} \Big\{ g \big( \frac{v}{\vt} \big) - \frac{v}{\vt} + 1 - \log \frac{\vt}{v_0} \Big\}(y,\tau) \, \dif y %\nonumber \\
\ge - C_0.
\end{align}
In deriving the first inequality in (\ref{f4}), we have also used the inequality (\ref{psi}). Substituting the estimates (\ref{f1})--(\ref{f4}) in (\ref{D1}) gives the desired estimate (\ref{f0}).
\end{proof}

%% estimate of A(x,t)
%%
%In order to obtain the pointwise positivity of $v(x,t)$,
%we derive the estimate of $A$ and $B$.
\begin{lemma}\label{lemma:VW}
Fix $(k,t)\in\mathbb{N}\times[0,T]$. Suppose that assumptions in Lemma \ref{lemma:v} hold. Then there exist positive constants $c_1=c_1(C_0)$ and $c_2=c_2(C)$ independent of $m$, $T$, $t$, and $k$ such that
\begin{subequations}
\begin{equation}\label{VW1}
e^{c_1(\tau_2-\tau_1)} \le \dfrac{V^1_{k,t}(\tau_2)}{V^1_{k,t}(\tau_1)} ,\qu
e^{-c_2|u_b|^2 (\tau_2-\tau_1)} \le \dfrac{W^1_{k,\,t}(x,\tau_2)}{W^1_{k,\,t}(x,\tau_1)} \le 1,
\end{equation}
for all $x\in \mI_t^{k-1} \cap [B(t),M_0-2] $ and $0\le \tau_1\le \tau_2 \le t$. Moreover,
\begin{equation}
 e^{c_1(\tau_2-\tau_1)} \le \dfrac{V^2_{k,t}(\tau_2)}{V^2_{k,t}(\tau_1)},\qu
1 \le \dfrac{W^2_{k,\,t}(x,\tau_2)}{W^2_{k,\,t}(x,\tau_1)} \le e^{c_2|u_b|^2 (\tau_2-\tau_1)},\label{VW2}
\end{equation}
for all $x\in\mJ_{t}^{k-1}\cap [M_0,M(t)]$ and $t_x\le \tau_1\le \tau_2 \le t$.
\end{subequations}
\end{lemma}
\begin{proof}
Applying Jensen's inequality and Lemma \ref{lemma:vInt} on (\ref{V1}), we obtain
\begin{equation*}
\log \frac{V^1_{k,\,t}(\tau_2)}{V^1_{k,\,t}(\tau_1)} = \frac{K \gamma}{\mu} \int_{\tau_1}^{\tau_2} \int_{\mI_{t}^k} v^{-\gamma} \, \dif y \dif s  \ge \frac{K \gamma}{\mu} \int_{\tau_1}^{\tau_2} \Big( \int_{\mI_{t}^k} v \, \dif y \Big)^{-\gamma} \dif s \ge \frac{K \gamma}{\mu} C_0^{-\gamma} (\tau_2 - \tau_1).
\end{equation*}
Setting $c_1\vcentcolon= K\gamma C_0^{-\gamma}/\mu$, we obtain the first inequality. Next, due to (\ref{stProp1}), we have the lower estimate for $W^1_{k,\,t}$ as
\begin{equation}\label{e2}
\log \frac{W^1_{k,\,t}(x,\tau_2)}{W^1_{k,\,t}(x,\tau_1)} = - \frac{(n-1)\gamma}{\mu} \int_{\tau_1}^{\tau_2} \!\! \int_x^{M_0} \zeta_{k,\,t} \frac{|\ut|^2}{r^n} \, \dif y \dif s \ge - \frac{3(n-1)\gamma}{2\mu} |u_b|^2 (\tau_2 - \tau_1).
\end{equation}
Let $c_2\vcentcolon=3(n-1)\gamma/(2\mu)$, we get the second inequality, which completes the proof.
\end{proof}

\begin{corollary}\label{corol:A}
Fix $(k,t)\in\mathbb{N}\times[0,T]$. Then there exists constants $c_1,\, c_2>0$ independent of $m$, $T$, $t$, $k$ such that for $x\in[B(t),M_0-2]\cap\mI_t^{k-1}$, and $0 \le \tau_1 \le \tau_2 \le t$,
\begin{equation*} 
e^{(c_1 - c_2 |u_b|^2 ) (\tau_2 - \tau_1)} \le \frac{A^1_{k,\,t}(x,\tau_2)}{A^1_{k,\,t}(x,\tau_1)}.
\end{equation*}
Moreover, for $x\in[M_0,M(t)]\cap\mJ_t^{k-1}$, and $t_x \le \tau_1 \le \tau_2 \le t$,
\begin{equation*} 
e^{ c_1 (\tau_2 - \tau_1) }  \le \frac{A^2_{k,\,t}(x,\tau_2)}{A^2_{k,\,t}(x,\tau_1)}.
\end{equation*}
\end{corollary}
%\begin{proof}
%Adding (\ref{e1}) to (\ref{e2}) yields
%\begin{align*}
%\log \frac{A_{k,\,t}(x,\tau_2)}{A_{k,\,t}(x,\tau_1)}& \ge \Big\{ \frac{K \gamma}{\mu} C_1^{\gamma} - \frac{3(n-1)\gamma}{2\mu} |u_b|^2 \Big\} (\tau_2 - \tau_1) %\nonumber \\
%=\vcentcolon (c_1 - c_2 |u_b|^2) (\tau_2 - \tau_1).
%\end{align*}
%\end{proof}

With the help of Lemmas \ref{lemma:D}--\ref{lemma:VW}, we obtain the point-wise estimate for $v(x,t)$ in $(x,t)\in\mL_{1}(T)\cup\mL_2(T)$. Once again, to avoid repetition, we only present the proof for $(x,t)\in \mL_1(T)$ since the case for $\mL_2(T)$ follows in a similar manner.

\begin{lemma}\label{lemma:vL12}
Suppose that assumptions in Lemma \ref{lemma:v} hold. Then there exists $\delta_{*}^0=\delta_{*}^0(\rho_0,u_0,\rho_+,\mu,\gamma,K,n)>0$ such that if $|u_b|\le \delta_{*}^0$, then there exists constants $0<\ul{v}_{*} < \ol{v}_{*} < \infty$ with $\ul{v}_{*}=\ul{v}_{*}(C_0)$ and $\ol{v}_*=\ol{v}_*(C_0)$, such that
\begin{equation*}
\ul{v}_{*} \le v(x,t) \le \ol{v}_{*}, \qu \text{for } \ (x,t)\in \mL_1(T)\cup\mL_2(T).
\end{equation*}
\end{lemma}
\begin{proof}
First, we obtain the upper bound of $v(x,t)$. Taking $|u_b|>0$ so small that $c_1 - c_2 |u_b| > c_1/2 \equiv c$ holds in Lemma \ref{lemma:VW}. Then,
\begin{equation}\label{g1}
\frac{A^1_{k,\,t}(x,\tau)}{A^1_{k,\,t}(x,t)} \le e^{-c(t - \tau)}, \quad A^1_{k,\,t}(x,t)^{-1}\le e^{-ct}
\end{equation}
for $x \in \mI_t^{k-1}\cap[B(t),M_0-2]$ and  $0 \le \tau \le t$. Substituting the estimates (\ref{f0}) and (\ref{g1}) in (\ref{form1}), we have
\begin{equation}\label{g3}
v(x,t)^\gamma \le \dfrac{C_0}{A^1_{k,\,t}(x,t)} + C_0 \int_0^t \frac{A^1_{k,\,t}(x,\tau)}{A^1_{k,\,t}(x,t)} \, \dif \tau \le C_0 e^{-ct} + C_0 \int_0^t e^{-c (t - \tau)} \, \dif \tau \le C_0.
\end{equation}
Since this estimate holds independent of
$k$, $t$, $m$, and $T$, we have the upper bound of $v$ uniformly in $(x,t)\in \mL_1(T)$. Namely, there exists a constant $\ol{v}_*=\ol{v}_*(C_0)>0$ such that $v(x,t) \le \ol{v}_*$ for an arbitrary $(x,t) \in \mL_1(T)$.

Next, we derive the lower bound of $v(x,t)$.
Let $x \in \mI_t^{k-1}\cap [B(t),M_0-2]$. We set $C_1\vcentcolon=c_1/3$ and $|u_b|^2 \le C_1/c_2$, where $c_1,\,c_2>0$ are constants from Lemma \ref{lemma:VW}, then
\begin{equation}\label{VWC1}
\dfrac{V_{k,\,t}^1(\tau)}{V_{k,\,t}^1(t)} \le e^{-3C_1(t-\tau)}, \qu 1 \le \dfrac{W_{k,\,t}^1(\tau)}{W_{k,\,t}^1(t)} \le e^{C_1(t-\tau)} \qu \text{for all } \ 0\le \tau \le t. 
\end{equation}
Applying Lemma \ref{lemma:D} on (\ref{form1}) and
using the boundedness of $v_0(x)$, we have
\begin{equation}\label{sqrtVW}
v(x,t)^{\gamma} \le \frac{C_0}{V^1_{k,\,t}(t) W^1_{k,\,t}(x,t)} + C_0 \int_0^t \sqrt{ \frac{V^1_{k,\,t}(\tau)}{V^1_{k,\,t}(t)} } \frac{W^1_{k,\,t}(x,\tau)}{W^1_{k,\,t}(x,t)} \cdot \sqrt{ \frac{V^1_{k,\,t}(\tau)}{V^1_{k,\,t}(t)} } \, \dif \tau. 
\end{equation}
Applying the Cauchy-Schwarz inequality on the second term of (\ref{sqrtVW}), and using the estimate (\ref{VWC1}), we have
\begin{equation}\label{h5}
v(x,t)^\gamma  \le C_0 e^{-2C_1 t} + C_0 \lambda (1 - e^{-C_1 t}) + \dfrac{C_0}{\lambda} \int_0^t \frac{V^1_{k,\,t}(\tau)}{V^1_{k,\,t}(t)} \, \dif \tau,
\end{equation}
where $\lambda\in(0,1]$ is an arbitrary constant to be determined later. Applying Jensen's inequality and Lemma \ref{lemma:vInt}, we see that
\begin{equation}\label{h6}
	\int_{\mI_t^{k-1}} v(y,t)^\gamma \, \dif y
	\ge \Big( \int_{\mI_t^{k-1}} v(y,t) \, \dif y
	\Big)^\gamma \ge C_{0}^{-\gamma} =\vcentcolon C_2 > 0,
\end{equation}
Therefore, integrating (\ref{h5}) with respect to $x \in \mI_t^{k-1}$, we have from (\ref{h6}) that for $t \in [0,T]$,
\begin{equation*}
\int_0^t \frac{V^1_{k,\,t}(\tau)}{V^1_{k,\,t}(t)} \, \dif \tau \ge \frac{\lambda}{C_0} \big( C_2 - C_0 \lambda - C_0 e^{-2C_1 t} + C_0 \lambda e^{-C_1 t} \big),
\end{equation*}
On the above inequality, we first fix a constant $\lambda \in (0,1]$. Then there exists $t_0=t_0(\lambda,C_0,C_1,C_2)>0$ which is independent of $k$, $t$, $m$, and $T$ such that
\begin{equation*}
C_0 \lambda + C_0 e^{-2C_1 t} \le \frac{C_2}{2}  \qu \text{for all} \ t\ge t_0.
\end{equation*}
Thus for all $t\ge t_0.$
\begin{equation}\label{vest}
\int_0^t \frac{V^1_{k,\,t}(\tau)}{V^1_{k,\,t}(t)} \, \dif \tau \ge \dfrac{\lambda}{C_0} \dfrac{C_2}{2} = \vcentcolon c_0.
\end{equation}
Thus, applying Lemmas \ref{lemma:D} on equality (\ref{form1}) in Lemma \ref{lemma:form1}, and using the estimate \eqref{VWC1}, it yields that for $t\ge t_0$,
\begin{align}\label{h8}
v(x,t)^{\gamma} & \ge C_0 \int_0^t  \frac{V^1_{k,\,t}(\tau)W^1_{k,\,t}(x,\tau)}{V^1_{k,\,t}(t)W^1_{k,\,t}(x,t)} \, \dif \tau \ge C_0 \int_0^t \frac{V^1_{k,\,t}(\tau)}{V^1_{k,\,t}(t)} \, \dif \tau \ge c_0 > 0.
\end{align}
To complete the proof, it suffices to show
the lower bound of $v$ for $ t \in [0, t_0]$. Multiplying $A^1_{k,\,t}(x,t) D^1_{k,\,t}(x,t)$ on (\ref{form1}), 
and 
then applying Lemma \ref{lemma:D} and (\ref{VWC1}), we obtain that
\begin{equation}\label{h9}
V^1_{k,\,t}(t) e^{-C_1 t} v(x,t)^\gamma \le C_0 v_0(x)^{\gamma} + C_0 \int_0^t V^1_{k,\,t}(\tau) \, \dif \tau.
\end{equation}
Integrating (\ref{h9}) in $x$ over $\mI_{t}^{k-1}$
and using (\ref{h6}), we have for $t\in[0,t_0]$,
\begin{equation}\label{h10}
V^1_{k,\,t}(t) e^{-C_1 t} \le C_0 + C_0 e^{ C_1 t_0 } \int_0^t V^1_{k,\,t}(\tau) e^{-C_1 \tau} \, \dif \tau,
\end{equation}
Applying Gr\"onwall's inequality on (\ref{h10}) over $t\in[0,t_0]$ yields that
\begin{equation}\label{h11}
V_{k,\,t}(t) \le C_0 \exp\big( C_1 t_0 + t_0 C_0 e^{C_1 t_0} \big) \le C_{t_0}, \qu \text{ for } \ t \in [0,t_0],
\end{equation}
where $C_{t_0}>0$ is a constant independent of $t$, $k$, $m$, and $T$. The combination of estimates (\ref{VWC1}) and (\ref{h11}) gives the upper bound $A_{k,\,t}(x,t) \le C_{t_0}.$
Therefore, substituting this inequality in (\ref{form1}) from Lemma \ref{lemma:form1}, and using Lemma \ref{lemma:D} gives
the lower bound of $v$ as
\begin{equation}\label{h13}
v(x,t)^\gamma \ge C_0^{-1} \dfrac{v_0(x)^{\gamma}}{A^1_{k,\,t}(x,t)} \ge C_0^{-1} C_{t_0}^{-1} > 0.
\end{equation}
This estimate is uniform with respect to $(x,t)\in\mL_1(T)$. Hence, the combination of (\ref{h8}) and (\ref{h13}) yields the uniform lower bound of $v(x,t)$ for $(x,t)\in\mL_1(T)$.
\end{proof}

Since $(M_0,t)\in \mL_2(T)$ for all $t\in[0,T]$, we have from Lemma \ref{lemma:vL12},
\begin{equation}
\ul{v}_*\le v(M_0,t)\le \ol{v}_* \qu \text{for all } \ t\in[0,T].
\end{equation}
Using $M_0$ as a reference point, the representation formula for the case $(x,t)\in \mL_3(T)$ can be obtained. This is stated in the following lemma, which can be derived using a similar proof to that of Lemma \ref{lemma:form1}:
\begin{lemma}\label{lemma:form3}
    Suppose assumptions in Lemma \ref{lemma:v} hold. Then
	\begin{equation*}
		v(x,t)^{\gamma} =
		\frac{v_0(x)^\gamma v_0(M_0)^{-\gamma} v(M_0,t)^{-\gamma} + \frac{K \gamma}{\mu} \int_0^t A_{3}(x,s) D_{3}(x,s) \dif s}{A_{3}(x,t) D_{3}(x,t)},
	\end{equation*}
	for all $(x,t)\in\mL_3(T)$, where
	\begin{alignat*}{2}
		A_{3}(x,t) &  :=
		\exp \Big( && \frac{K \gamma}{\mu}
		\int_0^{t} \!\! v(M_0,s)^{\gamma} \, \dif s
		- \frac{(n-1)\gamma}{\mu} \int_0^{t} \!\! \int_{x}^{M_0} \frac{|\ut|^2}{r^{n}} \, \dif y \dif s
		\Big), \\
		D_{3}(x,t) & :=
		\exp \Big( && \frac{\gamma}{\mu} \int_{x}^{M_0}
		\big\{ \frac{\psi}{r^{n-1}} - f(r) \big\}(y,s) \, \dif y\Big\vert_{s=0}^{s=t} + \frac{(n-1)\gamma}{\mu} \int_0^{t} \!\! \int_{x}^{M_0} \frac{\psi^2}{r^n} \, \dif y \dif s \Big).
	\end{alignat*}
\end{lemma}
We remark that $|M_0-x|\le 2$ for all $(x,t)\in \mL_3(T)$, thus there is no need to introduce ``cut-off'' functions with unit length as in the cases of $\mL_1(T)$ and $\mL_{2}(T)$.

Repeating the exact same arguments presented in Lemmas \ref{lemma:D}--\ref{lemma:vL12}, we can also obtain the upper and lower point-wise estimate for $v(x,t)$ for $(x,t)\in\mL_3(T)$:
\begin{lemma}\label{lemma:vL3}
Suppose that assumptions in Lemma \ref{lemma:v} hold. Then there exists $\delta_{**}^0=\delta_{**}^0(\rho_0,u_0,\rho_+,\mu,\gamma,K,n)>0$ such that if $|u_b|\le \delta_{**}^0$, then there exists constants $0<\ul{v}_{**} < \ol{v}_{**} < \infty$ with dependence $\ul{v}_{**}=\ul{v}_{**}(C_0)$ and $\ol{v}_{**}=\ol{v}_{**}(C_0)$ such that
\begin{equation*}
\ul{v}_{**}\le v(x,t) \le \ol{v}_{**}, \qu \text{for } \ (x,t)\in \mL_3(T).
\end{equation*}
\end{lemma}

If we set $\ul{v}\vcentcolon=\min\{\ul{v}_*,\,\ul{v}_{**}\}$, $\ol{v}\vcentcolon=\max\{\ol{v}_*,\,\ol{v}_{**}\}$, and $\delta_0\vcentcolon=\min\{\delta_{*}^0,\,\delta_{**}^0\}$, then Lemmas \ref{lemma:vL12} and \ref{lemma:vL3} imply that if $|u_b|\le \delta_0$ then \eqref{vulb} holds, and this completes the proof of Lemma \ref{lemma:v}.\hfill$\square$

Since the specific volume $v$ satisfies
\begin{equation}
	\underline{v} \le v(x,t) \le \bar{v}
	\quad \text{for} \
	(x,t) \in \mL(T),
	\label{bddv}
\end{equation}
one can verify with a simple calculation that $ G(v, \vt)$ is equivalent to $ \phi^2$. Namely,
$ C(\ul{v},\,\ol{v})^{-1} \phi^2 \le  G(v, \vt) \le  C(\ul{v},\,\ol{v}) \phi^2$. Therefore, the energy $\vE$ is equivalent to $\psi^2 + \phi^2$, that is,
\begin{equation}\label{L2vE}
C_0^{-1} (\psi^2 + \phi^2) \le \vE \le C_0 (\psi^2 + \phi^2).
\end{equation}
where $C_0>0$ is the constant specified in (\ref{constants}). As a consequence, we also have the following corollary from Lemma \ref{lemma:be} and Proposition \ref{prop:phiG}:
\begin{corollary}\label{corol:phiL2}
There exists $\delta_0=\delta(\rho_0,u_0,\rho_{+},\mu,\gamma,K,n)>0$ such that for $|u_b|\le \delta_0$,
\begin{equation*}
\sup\limits_{t\in[0,T]}\int_{S(t)}\!\! |\phi(x,t)|^2 \, \dif x + |u_b| \int_{0}^{T}\!\! |\phi(B(t),t)|^2\, \dif t + |u_b|^3 \iint_{\mathcal{L}(T)} \dfrac{\phi^2}{r^{3n-2}}\, \dif x \dif t \le C_0. 
\end{equation*}
\end{corollary}
\section{Higher order uniform estimates of \texorpdfstring{$(\phi,\psi)$}{(phi,psi)}}\label{sec:higher}\setcounter{equation}{0}
Using the point-wise estimate of $v$, we can derive the $H^1$ estimates of $(\phi,\psi)$ uniformly in $m\in\mathbb{N}$ and $T>0$, by exploiting the structure of \eqref{nsl}. The major difficulty in the proof is boundary integral terms which will unavoidably depend on $m\in\mathbb{N}$, and this is mainly due to the moving outer boundary $x=M(t)$. We avoid this problem by restricting the domain of integration with suitable cut-off functions. 
\subsection{Weighted \texorpdfstring{$H^1$}{H1} estimates for \texorpdfstring{$\phi$}{phi}}
In this subsection, the estimate of first order spatial derivative of $\phi=v-\vt$ is obtained using the special structure of the compressible Navier-Stokes equations in Lagrangian coordinate. First, we define:
\begin{equation}\label{F}
	F\vcentcolon= \mu \dfrac{\phi_x}{v} - \dfrac{\psi}{r^{n-1}}.
\end{equation}
To derive the evolution equation for $F$, we substitute \eqref{nsl1} in \eqref{nsl2}, then multiply both sides of the resulting equation by $r^{1-n}$ to get:
\begin{equation*}
\Big( \dfrac{u}{r^{n-1}} \Big)_t + (n-1) \dfrac{u^2}{r^n} - \gamma K v^{-\gamma} \dfrac{v_x}{v} = \mu \Big( \dfrac{v_x}{v} \Big)_t,
\end{equation*}
where we have used the identity $r_t = u$ from \eqref{Rdiff}. Next, we rewrite $v=\phi+\vt$, $u=\psi+\ut$ on the above equation, then by \eqref{st} and \eqref{Rdiff}, it follows that
\begin{subequations}\label{FtQ}
\begin{gather}
F_t + \dfrac{\gamma K}{\mu} \dfrac{F}{v^{\gamma}} = (n-1) \dfrac{\psi^2}{r^n} + \gamma\dfrac{p(v)-p(\vt)}{v-\vt} \dfrac{\partial_r \rt}{r^{n-1} \rt^2} \phi + Q \dfrac{\psi}{r^{n-1}}
\label{Ft} \\
\text{where } \ Q
:=
\partial_r \ut + (n-1) \dfrac{\ut}{r} - \dfrac{\gamma K}{\mu} v^{-\gamma} + \mu r^{n-1} \partial_r \big(\dfrac{\partial_r \rt}{r^{n-1}\rt^2}\big)
\label{Q}
\end{gather}
\end{subequations}
By Lemmas \ref{lemma:st}, \ref{lemma:v}, and Mean Value theorem, we have
\begin{equation}\label{Qest}
|Q|_\infty
\le C_0, \qu \big|\dfrac{p(v)-p(\vt)}{v-\vt}\big|_{\infty} \le C_0, \qu |\partial_r\rt(r)|\le C|u_b|^2 r^{-2n+1}.    
\end{equation}
The next lemma gives the estimate for $F$. Note that the spatial domain of integration is restricted in $x\in[B(t),M_0]$. This is done so to avoid boundary term from $x=M(t)$, which would inevitably depend on the domain parameter $m\in\mathbb{N}$. 
\begin{lemma}\label{lemma:phiH1}
Suppose $(v, u)$ is a solution to
\eqref{nsl}--\eqref{lcompa} in $\mL(T)$, then there exists $\delta_0=\delta_0(\rho_0,u_0,\rho_+,\mu,\gamma,K,n)>0$ such that if $|u_b|\le \delta_0$ then
\begin{equation*}
\sup\limits_{t\in[0,T]}\int_{B(t)}^{M_0}\!\! r^{2n-2} F^2(x,t) \, \dif x
+|u_b| \int_{0}^{T}\!\! \phi_x^2(B(t),t) \dif t  +\int_{0}^{T}\!\!\!\int_{B(t)}^{M_0} r^{2n-4} F^2 \, \dif x \dif t
\le C_0.
\end{equation*}
\end{lemma}
\begin{proof}
Multiplying (\ref{FtQ}) by $r^{2(n-2)}F$, we obtain
\begin{align*}
\big(\dfrac{r^{2(n-2)}F^2}{2}\big)_t + \dfrac{K\gamma}{\mu} \dfrac{r^{2(n-2)}F^2}{v^{\gamma}} =& (n-2) \dfrac{\psi+\ut}{r} r^{2(n-2)}F^2 + (n-1) \dfrac{\psi^2}{r^2} r^{n-2}F\\
&+ \gamma \dfrac{p(v)-p(\vt)}{v-\vt} \dfrac{\partial_r \rt }{r \rt^2} \phi r^{n-2} F + Q \dfrac{\psi}{r} r^{n-2} F. 
\end{align*}
Integrating the above in $x\in[B(\tau),M_0]$ then $\tau\in[0,t]$, choosing $|u_b|\le \frac{K\gamma}{10(n-2)\mu \ol{v}^{\gamma}}$, and using Lemmas \ref{lemma:v} and (\ref{Qest}), it follows that
\begin{align}\label{pH1Temp1}
&\int_{0}^{t}\!\!\! \int_{B(\tau)}^{M_0}\!\! \big(\dfrac{r^{2(n-2)}F^2}{2}\big)_{\tau}(x,\tau) \, \dif x \dif \tau + \dfrac{K\gamma}{\mu \ol{v}^{\gamma}} \int_{0}^{t}\!\!\!\int_{B(\tau)}^{M_0}\!\! r^{2(n-2)}F^2(x,\tau) \, \dif x \dif \tau \nonumber\\
%\le& \dfrac{K\gamma}{2\mu \ol{v}^{\gamma}}\int_{0}^{t}\!\!\!\int_{B(\tau)}^{M_0}\!\! r^{2(n-2)}F^2 \, \dif x \dif \tau + \dfrac{5(n-2)^2\mu \ol{v}^{\gamma}}{2K\gamma} \int_{0}^{t} \!\! \sup\limits_{y\in S(\tau)} \big|\dfrac{\psi}{r}\big|^2(y,\tau) \int_{B(\tau)}^{M_0}\!\! r^{2(n-2)} F^2 \,\dif x \dif \tau \nonumber\\ &+ \dfrac{5\gamma \mu \ol{v}^{\gamma}}{2\rt(1)^4 K} \big| \dfrac{p(v)-p(\vt)}{v-\vt} \big|_{\infty}^2 \int_{0}^{t} \!\!\! \int_{B(\tau)}^{M_0}\!\! C |u_b|^4 r^{-4n} \phi^2 \, \dif x \dif \tau \nonumber\\ & + \dfrac{5(n-1)^2\mu \ol{v}^{\gamma}}{2K\gamma} \int_{0}^{t} \!\! \sup\limits_{y\in S(\tau)} \big| \dfrac{\psi^2}{r^4} \big| \int_{B(\tau)}^{M_0}\!\! \psi^2 \, \dif x \dif \tau + \dfrac{5\mu \ol{v}^{\gamma}}{2 K \gamma \ul{v}} |Q|_{\infty}^2 \int_{0}^{t} \!\!\! \int_{B(\tau)}^{M_0} \!\! \dfrac{v\psi^2}{r^2} \, \dif x \dif \tau \nonumber\\
\le & \dfrac{K\gamma}{2\mu \ol{v}^{\gamma}}\int_{0}^{t}\!\! \int_{B(\tau)}^{M_0}\!\! r^{2(n-2)}F^2 \, \dif x \dif \tau + C_0 \int_{0}^{t} \!\! \sup\limits_{y\in S(\tau)} \dfrac{\psi^2}{r^2}(y,\tau) \int_{B(\tau)}^{M_0}\!\! r^{2(n-2)} F^2 \,\dif x \dif \tau  \nonumber\\
&+ C_0 \int_{0}^{t} \sup\limits_{y\in S(\tau)} \dfrac{\psi^2}{r^4}(y,\tau) \int_{B(\tau)}^{M_0}\!\! \psi^2 \,\dif x \dif \tau + C_0 \int_{0}^{t} \!\! \int_{B(\tau)}^{M_0}\!\! \Big\{ |u_b|^4 \dfrac{\phi^2}{r^{4n}} +  \dfrac{v \psi^2}{r^2} \Big\} \, \dif x \dif \tau .
\end{align}
Next, by Lemma \ref{lemma:be}, Corollary \ref{corol:infpsi}, and (\ref{L2vE}), we have from $n\ge 2$ that
\begin{subequations}\label{pH1Temp2}
\begin{align}
&\int_0^{t}\!\!\! \int_{B(\tau)}^{M_0} \!\! \Big\{ |u_b|^4 \dfrac{\phi^2}{r^{4n}} + \dfrac{v\psi^2}{r^2} \Big\}\, \dif x \dif \tau \le \!\! \int_0^{t} \!\!\! \int_{B(\tau)}^{M_0} \!\! \Big\{ C_0 |u_b|^4 \dfrac{G(v,\vt)}{r^{3n-2}} + \dfrac{v\psi^2}{r^2} \Big\}\, \dif x \dif \tau \le C_0,\\
&\int_{0}^{t} \sup\limits_{y\in S(\tau)} \dfrac{\psi^2}{r^4}(y,\tau) \int_{B(\tau)}^{M_0}\!\! \psi^2 \, \dif x \dif \tau \le C_0 \int_{0}^{t} \sup\limits_{y\in S(\tau)} \psi^2(y,\tau)\, \dif \tau \le C_0.
\end{align}
\end{subequations}
Furthermore, using Leibniz's integral rule, and boundary condition (\ref{lbc}), we have
\begin{align}\label{pH1Temp3}
&\int_0^{t}\!\!\!\int_{B(\tau)}^{M_0}\!\! \big( \dfrac{r^{2(n-2)}F^2}{2} \big)_{\tau}\, \dif x \dif \tau = \int_0^{t}\!\! \dfrac{\dif }{\dif \tau }\!\! \int_{B(\tau)}^{M_0}\!\! \dfrac{r^{2(n-2)}F^2}{2}(x,\tau)\, \dif x \dif \tau +  \int_{0}^{t}\!\! B^{\prime}(\tau) \dfrac{F^2}{2}(B(\tau),\tau)\, \dif \tau \nonumber \\
&= \int_{B(t)}^{M_0} \dfrac{r^{2(n-1)}F^2}{2}(x,t)\, \dif x - \int_{0}^{M_0}  \dfrac{r^{2(n-1)}F^2}{2}(x,0)\, \dif x + \mu^2 |u_b|\!\! \int_{0}^{t} \dfrac{\phi_x^2}{2v}(B(\tau),\tau) \dif \tau
\end{align}
Substituting (\ref{pH1Temp2})--(\ref{pH1Temp3}) in (\ref{pH1Temp1}), and using Lemma \ref{lemma:v} we have that
\begin{align}\label{pH1Temp4}
& \int_{B(t)}^{M_0}\! \dfrac{r^{2(n-2)}F^2}{2} (x,t) \, \dif x + \dfrac{\mu^2|u_b|}{2\ol{v}}  \!\! \int_{0}^{t}\!\! \phi_x^2 (B(\tau),\tau) \, \dif \tau + \dfrac{K\gamma}{2 \mu \ol{v}^{\gamma}} \int_{0}^{t}\!\!\!\int_{B(\tau)}^{M_0}\!\! r^{2(n-2)}F^2 \, \dif x \dif \tau \nonumber\\
&\le C_0 + C_0 \int_{0}^{t} \!\! \sup\limits_{y\in S(\tau)} \dfrac{\psi^2}{r^2}(y,\tau) \int_{B(\tau)}^{M_0}\!\! r^{2(n-2)} F^2 \,\dif x \dif \tau.
\end{align}
By Gr\"onwall's inequality, and Corollary \ref{corol:infpsi} it holds that for all $t\in[0,T]$
\begin{align*}
\int_{B(t)}^{M_0}\! \dfrac{r^{2(n-2)}F^2}{2} (x,t) \, \dif x \le C_0 \exp\Big( C_0 \int_0^{t} \sup\limits_{y\in S(\tau)} \dfrac{\psi^2}{r^2}(y,\tau)\, \dif \tau \Big) \le C_0 e^{C_0}
\end{align*}
Putting the above into (\ref{pH1Temp4}), and using Corollary \ref{corol:infpsi}, we obtain that for all $t\in[0,T]$,
\begin{align}\label{pH1Concl1}
\int_{B(t)}^{M_0}\! r^{2(n-2)}F^2 (x,t) \, \dif x + |u_b| \!\! \int_{0}^{t}\!\! \phi_x^2 (B(\tau),\tau) \, \dif \tau +  \int_{0}^{t}\!\!\!\int_{B(\tau)}^{M_0}\!\! r^{2(n-2)}F^2 \, \dif x \dif \tau \le C_0.
\end{align}

Next, we choose $|u_b|\le \frac{K\gamma}{10(n-1)\mu \ol{v}^{\gamma}}$. Multiplying (\ref{FtQ}) by $r^{2(n-1)}F$, integrating the resultant equality in $x\in [B(\tau),M_0]$, then by
Leibniz's integral rule, the Cauchy-Schwartz inequality, Lemmas \ref{lemma:be} and \ref{lemma:v}, we obtain that
\begin{subequations}\label{pH1Temp5}
\begin{align}
&\dfrac{\dif }{\dif \tau}\!\! \int_{B(\tau)}^{M_0}\! \dfrac{r^{2(n-1)}F^2}{2}\, \dif x + \dfrac{\mu^2|u_b|}{2\ol{v}} \phi_x^2(B(\tau),\tau) + \dfrac{K\gamma}{2\mu \ol{v}^{\gamma}}\!\! \int_{B(\tau)}^{M_0}\!\! r^{2(n-1)}F^2 \, \dif x
%&\le C_0 \int_{B(\tau)}^{M_0}\!\! \psi^2(x,\tau)\, \dif x + C_0 \int_{B(\tau)}^{M_0}\!\! r^{2n-4} \psi^2 F^2(x,\tau) \, \dif x + \mathcal{R}(\tau)\\
\le C_0 + \mathcal{D}(\tau)\label{pH1Temp5a}
\end{align}
\begin{align}
\text{where } \ \mathcal{D}(\tau) \vcentcolon=& C_0\!\! \sup\limits_{y\in S(\tau)} \dfrac{\psi^2}{r^2}(y,\tau) \int_{B(\tau)}^{M_0}\!\! \psi^2(x,\tau) \,\dif x + C_0 \!\! \int_{B(\tau)}^{M_0}\!\! |u_b|^4 \dfrac{\phi^2}{r^{4n-2}}(x,\tau) \, \dif x \nonumber\\
&+ C_0\!\! \int_{B(\tau)}^{M_0} r^{2n-4} \psi^2 F^2 (x,\tau)\, \dif x.\label{pH1Temp5b}
\end{align}
\end{subequations}
It follows from the previous estimate (\ref{pH1Concl1}) and Corollary \ref{corol:infpsi} that
\begin{equation}\label{pH1Temp6}
\int_{0}^{T}\!\!\int_{B(\tau)}^{M_0} r^{2n-4} \psi^2 F^2 \, \dif x \dif \tau \le \int_{0}^{T}\!\!\sup\limits_{y\in S(\tau)}\psi^2(y,\tau)  \int_{B(\tau)}^{M_0} r^{2(n-2)} F^2 \, \dif x \dif \tau \le C_0.
\end{equation}
Furthermore, by Lemma \ref{lemma:be}, Corollary \ref{corol:infpsi}, (\ref{L2vE}), and the fact that $n\ge 2$, we obtain
\begin{align}\label{pH1Temp7}
&\int_{0}^{T}\!\! \Big\{ \sup\limits_{y\in S(\tau)} \dfrac{\psi^2}{r^2}(y,\tau)\!\! \int_{B(\tau)}^{M_0} \!\! \psi^2(x,\tau) \, \dif x + \int_{B(\tau)}^{M_0} |u_b|^4 \dfrac{\phi^2}{r^{4n-2}}(x,\tau)\, \dif x \Big\} \dif \tau\nonumber\\
\le& C_0\int_{0}^{T} \sup\limits_{y\in S(\tau)} \psi^2(y,\tau) \dif \tau + C_0|u_b|^4 \int_{0}^{T}\!\!\int_{B(\tau)}^{M_0} \dfrac{G(v,\vt)}{r^{3n-2}}\, \dif x \dif \tau \le C_0.
\end{align}
Combining (\ref{pH1Temp6})--(\ref{pH1Temp7}), we conclude that $\int_{0}^{T} \mathcal{D}(\tau) \dif \tau \le C_0$. From the estimate (\ref{pH1Temp5a}), it follows that
\begin{equation}\label{pH1Temp8}
\dfrac{\dif}{\dif\tau } \mathcal{A}(\tau) + \dfrac{K\gamma}{\mu \ol{v}^{\gamma}} \mathcal{A}(\tau) \le C_0 + \mathcal{D}(\tau), \ \text{ where } \ \mathcal{A}(\tau)\vcentcolon = \int_{B(\tau)}^{M_0} \dfrac{r^{2(n-1)} F^2 }{2} (x,\tau) \, \dif x.
\end{equation}
To proceed, we use the following variant of the Gr\"onwall's inequality, whose proof is omitted as it is elemental.
\begin{proposition}\label{prop:Gron1}
Let $\mathcal{A}$, $\mathcal{B}$, $\mathcal{D}$ be non-negative functions such that $\mathcal{A}\in \mathcal{C}^1([0,T])$ and $\mathcal{B},\, \mathcal{D}\in L^1(0,T)$. Furthermore, let $\beta$, $C_0$ be two fixed positive constants. Suppose for all $\tau\in[0,T]$,
\begin{equation*}
\dfrac{\dif }{\dif \tau} \mathcal{A}(\tau) + \beta \mathcal{A}(\tau) \le C_0 + \mathcal{B}(\tau)\mathcal{A}(\tau) + \mathcal{D}(\tau) \qu \text{and } \int_{0}^{T}\big( \mathcal{B} +\mathcal{D} \big) (\tau)\, \dif \tau \le C_0.
\end{equation*}
Then $\mathcal{A}(\tau)\le e^{C_0} \big\{ \mathcal{A}(0) + C_0 (1 + \frac{1}{\beta}) \big\} $ for all $\tau\in[0,T]$.
\end{proposition}
Setting $\beta\vcentcolon= \frac{K\gamma}{\mu \ol{v}^{\gamma}}$ and $\mathcal{B}(t)=0$ in (\ref{pH1Temp8}) and applying Proposition \ref{prop:Gron1}, we obtain
\begin{align*}
&\sup\limits_{\tau\in[0,T]} \int_{B(\tau)}^{M_0} \dfrac{r^{2(n-1)} F^2 }{2} (x,\tau) \, \dif x \le \int_0^{M_0} \dfrac{r^{2(n-1)}F^2}{2}(x,0)\, \dif x + C_0 \big( 1 + \dfrac{\mu \ol{v}^{\gamma}}{K\gamma} \big)\\
&= \dfrac{1}{2}\int_0^{M_0} \big| \mu\dfrac{ r_0^{n-1} (\phi_0)_x}{v_0} - \psi_0 \big|^2 \, \dif x + C_0 \big( 1 + \dfrac{\mu \ol{v}^{\gamma}}{K\gamma} \big) \le C_0.
\end{align*}
The above estimate together with (\ref{pH1Concl1}) concludes the proof of this lemma.
\end{proof}

We recall that $F\vcentcolon= \mu v^{-1}\phi_x - r^{1-n}\psi$. Thus, combining Lemma \ref{lemma:phiH1} and the basic energy estimate, Lemma \ref{lemma:be}, it yields the following corollary.
\begin{corollary}\label{corol:phiH1}
Suppose $(v, u)$ is a solution to
\eqref{nsl}--\eqref{lcompa} in $\mL(T)$, then there exists $\delta_0=\delta_0(\rho_0,u_0,\rho_+,\mu,\gamma,K,n)>0$ such that if $|u_b|\le \delta_0$ then
\begin{equation}
\sup\limits_{t\in[0,T]}\int_{B(t)}^{M_0} r^{2n-2} \phi_x^2(x,t) \, \dif x +|u_b| \int_{0}^{T}\!\! \phi_x^2(B(t),t)\, \dif t
+ \int_{0}^{T}\!\!\!\int_{B(t)}^{M_0} r^{2n-4} \phi_x^2 \, \dif x \dif t
\le C_0. \label{i0}
\end{equation}
\end{corollary}

\subsection{Weighted \texorpdfstring{$H^1$}{H1} estimate of \texorpdfstring{$\psi$}{psi}}
We will derive the $L^2$-estimate of $\psi_x$ with the help of $H^1$-estimates of $\phi$ obtained in the previous subsection. Since the spatial domain of integration in Corollary \ref{corol:phiH1} is restricted in $x\in[B(t),M_0]$, it is necessary to introduce a suitable cut-off function. Therefore, we set $\tilde{\eta}\vcentcolon \R \to [0,1]$ as
\begin{equation}\label{teta}
\tilde{\eta}(y)=\left\{\begin{aligned}
&1 && \text{if $y\le -1$,}\\
&1-2(y+1)^2 && \text{if $-1\le y \le -1/2$,}\\
&2y^2 && \text{if $-1/2\le y \le 0$,}\\
&0 && \text{if $0\le y$.} 
\end{aligned}\right.
\end{equation}
One can verify that $\tilde{\eta}(y)\in\mathcal{C}^1(\R)$, and it satisfies:
\begin{equation}\label{tetaProp}
0\le \tilde{\eta}(y) \le 1, \qu |\tilde{\eta}^{\prime}(y)|^2 \le 8 \tilde{\eta}(y), \qu \supp(\tilde{\eta}^{\prime}) \subseteq [-1,0] \qu \text{for $y\in \R$.}
\end{equation}
Using $\tilde{\eta}$ and $R(x,t)$ from \eqref{R}, we set the cut-off function $\eta(x,t)\vcentcolon \mL(T)\to [0,1]$ as
\begin{equation}\label{eta}
	\eta(x,t) \vcentcolon= \tilde{\eta}(R(x,t)-R(M_0,t)), \qu \text{for $(x,t)\in\mL(T)$.}
\end{equation}
It then follows by \eqref{Rdiff} and \eqref{tetaProp} that for all $(x,t)\in \mL(T)$,
\begin{subequations}\label{etaProp}
\begin{align}
& 0 \le \eta(x,t)\le 1  \qu \text{and} \qu \eta(x,t) = \left\{ \begin{aligned}
&  1 && \text{ if } \ R(x,t)\le R(M_0,t)-1, \\
& 0 && \text{ if } \ R(M_0,t)\le R(x,t),
\end{aligned} \right. \\ 
&|\eta_x(x,t)| \le \dfrac{v\sqrt{8\eta}}{r^{n-1}}(x,t), \qu |\eta_t(x,t)| \le |u(x,t)-u(M_0,t)|\sqrt{8\eta(x,t)}
\end{align}
\end{subequations} 

Before proceeding to derive the $H^1$-estimate of $\psi$, we show the following proposition:
\begin{proposition}\label{prop:psiSob}
For each $\delta>0$ and $t\in[0,T]$,
\begin{subequations}\label{psiSob}
\begin{align}
\sup\limits_{x\in S(t)} \eta r^{3n-4} \psi_x^2 (x,t) \le C_0\dfrac{\delta+1}{\delta} \!\! \int_{S(t)} \!\!\! r^{2(n-1)} \psi_x^2 (x,t)\, \dif x + \delta\!\! \int_{S(t)}\!\!\! \eta  \dfrac{r^{4n-6} \psi_{xx}^2}{v} (x,t)\, \dif x,\label{psiSobA}\\
\sup\limits_{x\in S(t)} \eta r^{3n-3} \psi_x^2 (x,t) \le C_0\dfrac{\delta+1}{\delta}\!\! \int_{S(t)}\!\!\! r^{2(n-1)} \psi_x^2 (x,t)\, \dif x + \delta\!\! \int_{S(t)}\!\!\! \eta  \dfrac{r^{4n-4} \psi_{xx}^2}{v} (x,t)\, \dif x.\label{psiSobB}
\end{align}
\end{subequations}
\end{proposition}
\begin{proof}
Only the proof for (\ref{psiSobA}) is shown as (\ref{psiSobB}) follows in the exact same manner. Fix $t\in[0,T]$. By (\ref{etaProp}), it holds $\eta(M_0,t)=0$. Hence by Fundamental theorem of calculus and (\ref{Rdiff}),
we have for all $y\in S(t)$,
\begin{align}
&\eta r^{3n-4} \psi_x^2 (y,t) = \int_{M_0}^{y} \big( \eta r^{3n-4} \psi_x^2 \big)_x \dif x \nonumber\\
=& \int_{M_0}^{y} \big\{ \eta_x r^{3n-4} \psi_x^2 + (3n-4) \eta r^{2n-4} v \psi_x^2 + 2 \eta r^{3n-4} \psi_{xx} \psi_x \big\}(x,t) \dif x\nonumber
\end{align}
Using Lemma \ref{lemma:v}, (\ref{etaProp}), and 
the 
Cauchy-Schwarz inequality, it holds that
\begin{align*}
\eta r^{3n-4} \psi_x^2 (y,t) \le %&\big|\dfrac{v \sqrt{8\eta} }{r} \big|_{\infty} \int_{S(t)}  r^{2n-2} \psi_x^2 \, \dif x + (3n-4) \big| \dfrac{\eta v}{r^2} \big|_{\infty} \int_{S(t)}  r^{2n-2} \psi_x^2 \, \dif x\\ &+ \dfrac{\ol{v}}{\delta}\int_{S(t)} \eta r^{2n-2} \psi_x^2\, \dif x + \delta \int_{S(t)} \eta \dfrac{r^{4n-6} \psi_{xx}^2}{v} \, \dif x
\big( C_0 + \dfrac{\ol{v}}{\delta} \big) \int_{S(t)} r^{2n-2} \psi_x^2(x,t) \, \dif x + \delta \int_{S(t)} \eta \dfrac{r^{4n-6} \psi_{xx}^2}{v} \, \dif x.
\end{align*}
This proves (\ref{psiSobA}).
\end{proof}

\begin{lemma}\label{lemma:H1psi1}
For the solution $(v, u)$ to
\eqref{nsl}--\eqref{lcompa} in $\mL(T)$, it holds that
\begin{equation*}
\sup_{t\in[0,T]}\int_{S(t)} \eta r^{2n-4} \psi_x^2 (x,t) \, \dif x
+ \iint_{\mL(T)} \eta r^{4n-6} \psi_{xx}^2 \, \dif x \dif t
\le C_0.
\end{equation*}
\end{lemma}
\begin{proof}
Using (\ref{Rdiff}), we compute to get that
\begin{align}\label{uH1Temp1}
&-\int_{S(t)} \eta r^{2n-4} \psi_t \psi_{xx} \dif x  \nonumber\\
=& \int_{S(t)}\!\! \Big\{ \Big( \eta r^{2n-4} \dfrac{\psi_x^2}{2} \Big)_t -\big( \eta r^{2n-4} \psi_t \psi_x \big)_x - r^{2n-4} \eta_t \dfrac{\psi_x^2}{2} \Big\}\,  \dif x  \nonumber\\
&+\int_{S(t)}\!\! \big\{  (2n-4) \eta v r^{n-4} \psi_x \psi_t + \eta_x r^{2n-4} \psi_x \psi_t - (n-2) \eta r^{2n-5} u \psi_x^2 \big\}\, \dif x.
\end{align}
By (\ref{B}) and the boundary condition (\ref{lbc}) it follows that 
\begin{equation*}
\psi_t(B(t),t) + |u_b| \dfrac{\psi_x}{v}(B(t),t) = 0 \qu \text{for all } \ t\in[0,T].
\end{equation*}
Thus using the above equality, (\ref{etaProp}), and Leibniz's integral rule, we have that
\begin{align*}
&-\int_{S(t)}\!\!\big( \eta r^{2n-4} \psi_t \psi_x \big)_x \dif x %= (\eta \psi_\tau \psi_x) (B(t),t) 
= - |u_b| \dfrac{ \eta \psi_x^2}{v} (B(\tau),\tau),\nonumber\\ 
&\int_{S(t)}\!\! \Big( \eta r^{2n-4} \dfrac{\psi_x^2}{2} \Big)_t (x,t)\, \dif x = \dfrac{\dif }{\dif t} \int_{S(t)} \eta\dfrac{r^{2n-4}\psi_x^2}{2}(x,t) \, \dif x + \dfrac{|u_b|}{2} \dfrac{\eta \psi_x^2}{v}(B(t),t)
\end{align*}
Substituting the above identities in (\ref{uH1Temp1}), we have
\begin{align}\label{uH1Temp2}
&-\int_{S(t)} \eta r^{2n-4} \psi_t \psi_{xx} \dif x  \nonumber\\
=& \dfrac{\dif }{\dif t}\! \int_{S(t)}\!\! \eta \dfrac{r^{2n-4}\psi_x^2}{2} (x,t) \dif x  - \dfrac{|u_b|}{2} \dfrac{\eta \psi_x^2}{v}(B(t),t)  -\int_{S(t)}\!\!  r^{2n-4} \eta_t \dfrac{\psi_x^2}{2} \,  \dif x  \nonumber\\
&+\int_{S(t)}\!\! \big\{  (2n-4) \eta v r^{n-4} \psi_x \psi_t + \eta_x r^{2n-4} \psi_x \psi_t - (n-2) \eta r^{2n-5} u \psi_x^2 \big\}\, \dif x.
\end{align}

Next, rewriting (\ref{nsl2}) using (\ref{Rdiff}), we have 
\begin{align}\label{psieq}
&\begin{aligned}
\psi_t - \mu \dfrac{r^{2(n-1)}\psi_{xx}}{v} =& 2(n-1) \mu r^{n-2} \psi_x - \mu \dfrac{r^{2(n-1)}}{v^2} \phi_x \psi_x + \mu \dfrac{\partial_r \rt }{\rt^2} \dfrac{r^{n-1}}{v} \psi_x\\
&+ \gamma p(v) \dfrac{r^{n-1}}{v} \phi_x - (n-1)\mu \dfrac{v \psi}{r^2}  + \mathcal{R},\\
\end{aligned}\\
& \text{where} \qu \mathcal{R} = \dfrac{\partial_r \rt }{\rt^2 }  \dfrac{p^{\prime}(v)-p^{\prime}(\vt)}{v-\vt} v \phi - \partial_r \ut \psi  + \rt \ut \partial_r \ut \phi.  \nonumber
\end{align}
Multiplying (\ref{psieq}) by $-\eta r^{2n-4} \psi_{xx} $, integrating the resultant equality in $x\in S(t)$, then substituting the identity (\ref{uH1Temp2}), we obtain that
\begin{align}\label{uH1Temp3}
& \dfrac{\dif }{\dif t}\! \int_{S(t)}\!\! \eta \dfrac{r^{2n-4}\psi_x^2}{2} (x,t) \dif x + \mu \int_{S(t)} \eta \dfrac{r^{4n-6}\psi_{xx}^2}{v} (x,t)\, \dif x \nonumber\\
%= & \dfrac{|u_b|}{2} \dfrac{\eta \psi_x^2 }{v}(B(t),t) +\!\! \int_{S(t)}\!\!\! \big\{ (n-2) \eta r^{2n-5} u \psi_x^2 - (2n-4) \eta v r^{n-4} \psi_x \psi_t - \eta_x r^{2n-4} \psi_x \psi_t  \big\}\, \dif x \nonumber\\ &+\int_{S(t)}\!\! \Big\{ r^{2n-4} \eta_t \dfrac{\psi_x^2}{2} - 2(n-1)\mu \eta r^{3n-6} \psi_x \psi_{xx} + \mu \eta \dfrac{r^{4n-6}}{v^2}\phi_x \psi_x \psi_{xx} - \eta r^{2n-4} \psi_{xx} \mathcal{R}  \Big\} \,  \dif x \nonumber\\ &+\!\! \int_{S(t)}\!\! \Big\{ (n-1) \mu \eta r^{2n-6} v \psi_{xx} \psi - \mu \eta \dfrac{\partial_r \rt}{\rt^2} \dfrac{r^{3n-5}}{v} \psi_x \psi_{xx} - \gamma \eta \dfrac{p(v)}{v} r^{3n-5} \phi_x \psi_{xx}  \Big\} \, \dif x \nonumber\\
=& \dfrac{|u_b|}{2} \dfrac{\eta \psi_x^2 }{v}(B(t),t) - \!\! \int_{S(t)}\!\!\! \big\{ (2n-4) \eta v r^{n-4} \psi_x \psi_t + \eta_x r^{2n-4} \psi_x \psi_t  \big\}\, \dif x \nonumber\\
&+\int_{S(t)} \Big\{ r^{2n-4} \eta_t \dfrac{\psi_x^2}{2}  + (n-2) \eta r^{2n-5} u \psi_x^2 \Big\}\, \dif x + \mu \int_{S(t)} \eta \dfrac{r^{4n-6}}{v^2} \phi_x \psi_x \psi_{xx} \, \dif x \nonumber\\
&- \gamma \int_{S(t)} \eta \dfrac{p(v)}{v} r^{3n-5} \phi_x \psi_{xx} \, \dif x - \mu \int_{S(t)} \Big\{ \dfrac{2(n-1)}{r} + \dfrac{\partial_r \rt}{v \rt^2}  \Big\} \eta r^{3n-5} \psi_x \psi_{xx}\, \dif x \nonumber\\
&+\int_{S(t)} \Big\{ (n-1)\mu\dfrac{v \psi}{r^2} - \mathcal{R} \Big\} \eta r^{2n-4} \psi_{xx} \, \dif x =\vcentcolon \sum\limits_{i=1}^{7} \mathcal{U}_i
%\text{\textbf{\textcolor{red}{not finished!}}}
\end{align}
Using Proposition \ref{prop:psiSob}, we estimate $\mathcal{U}_1$ as:
\begin{align*}
\mathcal{U}_1\vcentcolon= \dfrac{|u_b|}{2} \dfrac{\eta \psi_x^2 }{v}(B(t),t) \le C_0 \int_{S(t)} \!\!\! r^{2(n-1)} \psi_x^2 \, \dif x + \dfrac{\mu}{12} \int_{S(t)}\!\!\! \eta \dfrac{r^{4n-6}\psi_{xx}^2}{v} \, \dif x
\end{align*}
Since $0\le \eta \le 1$, it holds that $\eta\le \sqrt{\eta}$. Thus using (\ref{stProp3}), then substituting (\ref{psieq}), we have by Lemma \ref{lemma:v} and
the 
Cauchy-Schwarz inequality:
\begin{align}\label{uH1Temp4}
&\mathcal{U}_2 \vcentcolon= - \int_{S(t)}\!\!\! \big\{ (2n-4) \eta v r^{n-4} \psi_x \psi_t + \eta_x r^{2n-4} \psi_x \psi_t  \big\}\, \dif x %\nonumber\\ \le& \int_{S(t)}\!\!\! \big\{ (2n-4) \dfrac{\eta}{r} + \sqrt{8\eta} \big\} v r^{n-3}|\psi_x \psi_t|\, \dif x \nonumber\\
\le C \int_{S(t)}\!\!\! \sqrt{\eta} v r^{n-3} |\psi_x \psi_t| \, \dif x\nonumber\\
%\le & C_0 \int_{S(t)}\!\!\! \sqrt{\eta} r^{n-3} |\psi_x| \Big\{ r^{2n-2} |\psi_{xx}| + r^{n-2} |\psi_x| + r^{2n-2} \phi_x \psi_x + r^{-n} \psi_x + r^{n-1}\phi_x + \dfrac{\psi}{r^2} + \mathcal{R}  \Big\} \, \dif x \nonumber\\
\le & C_0 \int_{S(t)}\!\!\! \sqrt{\eta} r^{n-3} |\psi_x| \Big\{ r^{2n-2} |\psi_{xx}| + r^{n-2} |\psi_x| + r^{2n-2} |\phi_x \psi_x|  + r^{n-1} |\phi_x| + \dfrac{|\psi|}{r^2} + |\mathcal{R}|  \Big\} \, \dif x \nonumber\\
\le & \dfrac{\mu}{24} \int_{S(t)}\!\!\! \eta \dfrac{r^{4n-6}\psi_{xx}^2}{v} \, \dif x + \int_{S(t)} \eta r^{4n-6}\psi_x^2 \phi_x^2\, \dif x  + C_0\int_{S(t)}\!\!\! \big\{ r^{2n-2} \psi_x^2 + \dfrac{v \psi^2}{r^2} \big\} \, \dif x\nonumber\\
&+ C_0 \int_{S(t)} \eta r^{2n-6} \phi_x^2 \, \dif x  + C_0 \int_{S(t)} \mathcal{R}^2 \, \dif x.
\end{align}
By (\ref{etaProp}), $\supp\big(\eta(\cdot,t)\big)\subseteq [B(t),M_0] $ for each $t\in[0,T]$. Using Lemma \ref{lemma:phiH1} and Proposition \ref{prop:psiSob}, it follows that
\begin{align}\label{uH1Temp5}
&C_0 \int_{S(t)}\!\!\! \eta r^{4n-6} \psi_x^2 \phi_x^2\, \dif x  \le C_0 \sup_{y\in S(t)} \eta r^{3n-4} \psi_x^2 (y,t) \int_{B(t)}^{M_0}\!\!\! r^{n-2} \phi_x^2 \, \dif x\nonumber\\ 
&\le C_0  \sup_{y\in S(t)} \eta r^{3n-4} \psi_x^2 (y,t) \le \dfrac{\mu}{24} \int_{S(t)}\!\!\! \eta \dfrac{r^{4n-6}\psi_{xx}^2}{v} \, \dif x + C_0 \int_{S(t)}\!\!\! r^{2n-2} \psi_x^2\, \dif x. 
\end{align}
Thus substituting (\ref{uH1Temp5}) in (\ref{uH1Temp4}), we obtain that
\begin{align*}
\mathcal{U}_2 \le & \dfrac{\mu}{12} \int_{S(t)}\!\!\! \eta \dfrac{r^{4n-6}\psi_{xx}^2}{v} \, \dif x + C_0 \int_{B(t)}^{M_0}\!\! r^{2n-4} \phi_x^2 \, \dif x + C_0\int_{S(t)}\!\!\! \big\{ r^{2n-2} \psi_x^2 + \dfrac{v \psi^2}{r^2} + \mathcal{R} \big\} \, \dif x.
\end{align*}
For $\mathcal{U}_3$, we use (\ref{stProp1}), (\ref{etaProp}), and $\eta\le \sqrt{\eta}$ to get
\begin{align*}
\mathcal{U}_3 \vcentcolon=& \int_{S(t)}\!\!\! \Big\{ \dfrac{\eta_t}{2} r^{2n-4} \psi_x^2 + (n-2) \eta r^{2n-5} u \psi_x^2 \Big\}\, \dif x \le C \sup_{y\in S(t)} |u(y,t)| \int_{S(t)} \sqrt{\eta} r^{2n-4} \psi_x^2\, \dif x\\
\le & C \sup_{y\in S(t)} |\psi(y,t)| \int_{S(t)} \sqrt{\eta} r^{2n-4} \psi_x^2 \, \dif x + C \sup_{y\in S(t)} |\ut(y,t)| \int_{S(t)} \sqrt{\eta} r^{2n-4} \psi_x^2 \, \dif x\\
\le & C \sup_{y\in S(t)} \psi^2(y,t) \int_{S(t)} \eta r^{2n-4} \psi_x^2 \, \dif x + C \int_{S(t)} r^{2n-2} \psi_x^2\, \dif x. 
\end{align*}
The term $\mathcal{U}_4$ is estimated by Lemma \ref{lemma:v}, 
the 
Cauchy-Schwarz inequality, and (\ref{uH1Temp5})
\begin{align*}
&\mathcal{U}_4 \vcentcolon= \mu \int_{S(t)}\!\!\! \eta \dfrac{r^{4n-6}}{v^2} \phi_x \psi_x \psi_{xx} \, \dif x \le \dfrac{\mu}{24}\int_{S(t)}\!\!\! \eta \dfrac{r^{4n-6}\psi_{xx}^2}{v} \, \dif x + C\int_{S(t)}\!\!\! \eta r^{4n-6} \psi_x^2 \phi_x^2\, \dif x\\
&\le \dfrac{\mu}{12}\int_{S(t)}\!\!\! \eta \dfrac{r^{4n-6}\psi_{xx}^2}{v} \, \dif x + C_0 \int_{S(t)} \!\!\! r^{2n-2} \psi_x^2 \, \dif x.   
\end{align*}
To estimate the term $\mathcal{U}_5$, we use 
the 
Cauchy-Schwarz inequality, Lemma \ref{lemma:v}, and the fact that $\supp\big( \eta(\cdot,t) \big)\subseteq [B(t),M_0]$, to obtain that
\begin{align*}
\mathcal{U}_5 \le \dfrac{\mu}{12} \int_{S(t)} \eta \dfrac{r^{4n-6}\psi_{xx}^2}{v} \,\dif x + C_0 \int_{B(t)}^{M_0} r^{2n-4} \phi_x^2 \, \dif x.
\end{align*}
Using 
the 
Cauchy-Schwartz inequality, (\ref{stProp3}) and Lemma \ref{lemma:v}, the terms $\mathcal{U}_6$ and $\mathcal{7}$ are estimated as:
\begin{align*}
|\mathcal{U}_6|+|\mathcal{U}_7| \le \dfrac{\mu}{6} \int_{S(t)} \eta \dfrac{r^{4n-6}\psi_{xx}^2}{v} \,\dif x + C_0 \int_{S(t)} \Big\{ r^{2n-2} \psi_x^2 + \dfrac{v\psi^2}{r^2} + \mathcal{R}^2  \Big\}\, \dif x.
\end{align*}
Substituting the estimates for $\mathcal{U}_1$-$\mathcal{U}_7$ in (\ref{uH1Temp3}), it follows that
\begin{multline}
\label{uH1Temp6}
\dfrac{\dif }{\dif t}\! \int_{S(t)}\!\! \eta \dfrac{r^{2n-4}\psi_x^2}{2} (x,t) \dif x + \dfrac{\mu}{2} \int_{S(t)} \eta \dfrac{r^{4n-6}\psi_{xx}^2}{v} (x,t)\, \dif x \\
\le C \sup_{y\in S(t)} \psi^2(y,t) \int_{S(t)} \eta r^{2n-4} \psi_x^2 (x,t) \, \dif x + C_0 \int_{B(t)}^{M_0} r^{2n-4}\phi_x^2(x,t)\,\dif x \\
+ C_0 \int_{S(t)} \Big\{ r^{2n-2} \psi_x^2 + \dfrac{v\psi^2}{r^2} + \mathcal{R}^2  \Big\}(x,t)\, \dif x
\end{multline}
From the definition (\ref{psieq}), it holds that $\int_{0}^{T}\int_{S(t)} \mathcal{R}^2 \,\dif x \dif t \le C_0$ by (\ref{stProp3}) and Lemmas \ref{lemma:be}, \ref{lemma:v}. Therefore integrating (\ref{uH1Temp6}) in time, and using Lemmas \ref{lemma:be}, \ref{lemma:phiH1}, we have that for $t\in[0,T]$,

\begin{multline*}
\int_{S(t)}\!\! \eta \dfrac{r^{2n-4}\psi_x^2}{2} (x,t) \dif x + \dfrac{\mu}{2} \int_{0}^{t}\int_{S(\tau)} \eta \dfrac{r^{4n-6}\psi_{xx}^2}{v} \, \dif x \dif \tau\\
\le C_0 + C \int_{0}^{t}\!\! \sup_{y\in S(\tau)} \psi^2(y,\tau) \int_{S(\tau)} \eta r^{2n-4} \psi_x^2 \, \dif x \dif \tau.
\end{multline*}
Thus by Gr\"onwall's inequality and Corollary \ref{corol:infpsi}, we conclude that for all $t\in[0,T]$,
\begin{align*}
\int_{S(t)}\!\! \eta \dfrac{r^{2n-4}\psi_x^2}{2} (x,t) \dif x + \dfrac{\mu}{2} \int_{0}^{t}\int_{S(\tau)} \eta \dfrac{r^{4n-6}\psi_{xx}^2}{v} \, \dif x \dif \tau \le C_0 + C_0 e^{C_0}.
\end{align*}
Taking supremum in $t\in[0,T]$ on the above 
inequality yields the desired result.
\end{proof}
%\bigskip

\begin{lemma}
\label{lemma:H1psi2}
For the solution $(v, u)$ to
\eqref{nsl}--\eqref{lcompa} in $\mL(T)$, it holds that
\begin{equation*}
\sup\limits_{t\in[0,T]}\int_{S(t)} \eta r^{2n-2} \psi_x^2(x,t) \, \dif x
\le C_0.
\end{equation*}
\end{lemma}
\begin{proof}
Multiplying (\ref{psieq}) by $-\eta r^{2n-2} \psi_{xx} $, integrating the resultant equality in $x\in S(t)$, and repeating the same calculations as (\ref{uH1Temp1})--(\ref{uH1Temp3}) in the proof of Lemma \ref{lemma:H1psi1}, we obtain that
\begin{align}\label{psiH1Temp3}
& \dfrac{\dif }{\dif t}\! \int_{S(t)}\!\! \eta \dfrac{r^{2n-2}\psi_x^2}{2} (x,t) \dif x + \mu \int_{S(t)} \eta \dfrac{r^{4n-4}\psi_{xx}^2}{v} (x,t)\, \dif x \nonumber\\
=& \dfrac{|u_b|}{2} \dfrac{\eta \psi_x^2 }{v}(B(t),t) - \!\! \int_{S(t)}\!\!\! \big\{ 2(n-1) \eta v r^{n-2} \psi_x \psi_t + \eta_x r^{2n-2} \psi_x \psi_t  \big\}\, \dif x \nonumber\\
&+\int_{S(t)} \Big\{ r^{2n-2} \eta_t \dfrac{\psi_x^2}{2}  + (n-1) \eta r^{2n-3} u \psi_x^2 \Big\}\, \dif x + \mu \int_{S(t)} \eta \dfrac{r^{4n-4}}{v^2} \phi_x \psi_x \psi_{xx} \, \dif x \nonumber\\
&- \gamma \int_{S(t)} \eta \dfrac{p(v)}{v} r^{3n-3} \phi_x \psi_{xx} \, \dif x - \mu \int_{S(t)} \Big\{ \dfrac{2(n-1)}{r} + \dfrac{\partial_r \rt}{v \rt^2}  \Big\} \eta r^{3n-3} \psi_x \psi_{xx}\, \dif x \nonumber\\
&+\int_{S(t)} \Big\{ (n-1)\mu\dfrac{v \psi}{r^2} - \mathcal{R} \Big\} \eta r^{2n-2} \psi_{xx} \, \dif x.
\end{align}
The right hand side of (\ref{psiH1Temp3}) has the same terms as that of (\ref{uH1Temp3}), except the exponent of $r$ in each integrand of (\ref{psiH1Temp3}) is less by a factor $r^{-2}$. Thus one can apply similar estimates as $\mathcal{U}_1$-$\mathcal{U}_7$ in the proof of Lemma \ref{lemma:H1psi1}, to obtain that
\begin{align}\label{psiH1Temp4}
&\dfrac{\dif }{\dif t} \int_{S(t)}\!\!\! \eta \dfrac{r^{2n-2}\psi_x^2}{2}(x,t)\, \dif x + \dfrac{\mu}{2}\int_{S(t)} \eta \dfrac{r^{4n-4}\psi_{xx}^2}{v}(x,t)\, \dif x,\nonumber\\
\le& C \sup_{y\in S(t)} \psi^2(y,t) \int_{S(t)}\!\!\! \eta r^{2n-2} \psi_x^2 (x,t)\, \dif x + C_0 \int_{B(t)}^{M_0} r^{2n-2} \phi_x^2(x,t)\, \dif x\nonumber\\
&+ C_0 \int_{S(t)} \Big\{ r^{2n-2}\psi_x^2 + \dfrac{v\psi^2}{r^2} + \mathcal{R}^2 \Big\}\, \dif x \nonumber\\
\le & C_0 + C\sup_{y\in S(t)} \psi^2(y,t) \int_{S(t)}\!\!\! \eta r^{2n-2} \psi_x^2 \, \dif x + C_0 \int_{S(t)}\!\!\! \Big\{ r^{2n-2}\psi_x^2 + \dfrac{v\psi^2}{r^2} + \mathcal{R}^2 \Big\}\, \dif x, 
\end{align}
where in the last line we used Lemma \ref{lemma:phiH1}. Next, using (\ref{beConcl}) from the proof of Lemma \ref{lemma:be}, Lemma \ref{lemma:v}, and the fact that $\eta \le 1$, it follows that
\begin{align}\label{psiH1Temp5}
\dfrac{\dif }{\dif t}\int_{S(t)}\!\! \vE(x,t)\, \dif x +  \int_{S(t)} \!\! \big\{  \dfrac{\mu}{\ol{v}} \eta r^{2n-2} \psi_x^2 + \vE \big\} (x,t)\, \dif x 
\le \int_{S(t)}\vE(x,t) \,\dif x \le C_0.
\end{align}
Multiplying (\ref{psiH1Temp4}) by $\frac{2\mu}{\ol{v}}$, then adding it to (\ref{psiH1Temp5}), we obtain
\begin{align*}
&\dfrac{\dif }{ \dif t} \mathcal{A}(t) + \mathcal{A}(t) \le C_0 + \mathcal{B}(t) \mathcal{A}(t) + \mathcal{D}(t), \\
\text{where} \quad & \mathcal{A}(t)\vcentcolon= \int_{S(t)} \big\{ \dfrac{\mu}{\ol{v}} \eta r^{2n-2} \psi_x^2 + \vE \big\}(x,t) \, \dif x, \qu \mathcal{B}(t)\vcentcolon= \sup_{y\in S(t)}\psi^2(y,t),\\
&\mathcal{D}(t)\vcentcolon = C_0 \int_{S(t)} \big\{ r^{2n-2}\psi_x^2 + \dfrac{v\psi^2}{r^2} + \mathcal{R}^2 \big\}\, \dif x. 
\end{align*}
By Lemma \ref{lemma:be} and Corollary \ref{corol:infpsi}, it follows that $\int_{0}^{T}(\mathcal{B}+\mathcal{D})(t)\dif t\le C_0.$ Thus applying Proposition \ref{prop:Gron1}, we have that
\begin{align*}
\int_{S(t)}\!\!\! \big\{ \dfrac{\mu}{\ol{v}}\eta r^{2n-2} \psi_x^2 + \vE \big\}(x,t)\, \dif x =\mathcal{A}(t) \le e^{C_0} \big\{ \mathcal{A}(0) + 2C_0 \big\} \le C_0,
\end{align*}
for all $t\in[0,T]$, which proves the lemma.
\end{proof}

Lemmas \ref{lemma:H1psi1}--\ref{lemma:H1psi2} give the uniform bound of $\psi(x,t)$.
\begin{corollary}
\label{corol:psibounds}
For the solution $(v, u)$ to \eqref{nsl}--\eqref{lcompa} in $\mL(T)$, it holds that
\begin{equation}
\sup_{(x,t)\in \mL(T)} | \eta r^{n-1} \psi^2|(x,t)
\le C_0.
\end{equation}
\end{corollary}
\begin{proof}
The computation similar to the derivation of (\ref{b3}) gives
\begin{align*}
| \eta r^{n-1} \psi^2 |
& = \Big| \int_{B(t)}^x ( \eta r^{n-1} \psi^2)_x \, dx \Big| \\
&\le \int_0^\infty \left| \eta_x r^{n-1} \psi^2 + (n-1) \eta \frac{v}{r} \psi^2  + 2 \eta r^{n-1} \psi \psi_x \right| \, \dif x \\
& \le C \int_{S(t)} \psi^2 \, \dif x + C \int_{S(t)} \eta r^{2n-2} |\psi_x|^2 \, \dif x
\le C_0.
\end{align*}
Here, the last inequality holds owing to Lemmas \ref{lemma:be} and \ref{lemma:H1psi2}.
\end{proof}

\begin{lemma}\label{lemma:Dt}
For the solution $(v, u)$ to \eqref{nsl}--\eqref{lcompa} in $\mL(T)$, it holds that
\begin{equation}
\iint_{\mL(T)} \phi_t^2 \, \dif x \dif t
\le C_0, \qu \iint_{\mL(T)} \eta|\psi_t|^2 \, \dif x \dif t
\le C_0 (1+T). \label{Dt}
\end{equation}
\end{lemma}
\begin{proof}
Rewriting equation (\ref{nsl1}) with $v=\phi + \vt$, using (\ref{st}) and (\ref{Rdiff}), we get:
\begin{equation*}
\phi_t = \dfrac{\partial_r \rt}{\rt^2} \psi + r^{n-1} \psi_x + (n-1)\dfrac{v\psi}{r} + \dfrac{\partial_r (r^{n-1}\ut)}{r^{n-1}}\phi.
\end{equation*}
Taking $\lVert \cdot \rVert_{L^2(\mL(T))}^2$ on both sides of the above equation, and using (\ref{stProp2})--(\ref{stProp3}) and Lemma \ref{lemma:be}, we obtain the first estimate of (\ref{Dt}). 

Next, multiplying $\eta$ on (\ref{psieq}) and taking $\lVert \cdot \rVert_{L^2(\mL(T))}^2$ on both sides of the resultant equation, we have by Lemmas \ref{lemma:be}, \ref{lemma:v}, \ref{lemma:phiH1}, and \ref{lemma:H1psi2} that
\begin{align*}
\iint_{\mL(T)}\eta |\psi_t|^2\, \dif x \dif t \le C_0(1+T) + C_0 \int_{0}^{T} \sup_{y\in S(t)} \eta r^{2n-2} \psi_x^2 (y,t) \, \dif t
\end{align*}
By Proposition \ref{prop:psiSob}, and Lemma \ref{lemma:H1psi2}, it follows that
\begin{align*}
\iint_{\mL(T)}\!\!\!\eta |\psi_t|^2\, \dif x \dif t \le& C_0(1+T) + C_0\iint_{\mL(T)}\!\!\! \Big\{ r^{2n-2} \psi_x^2 + \eta \dfrac{r^{4n-4}\psi_{xx}^2}{v} \Big\} \, \dif x \dif t \le C_0(1+T),
\end{align*}
which is the second estimate of (\ref{Dt}).
\end{proof}

\subsection{The H\"{o}lder estimate of \texorpdfstring{$(\rho,u)$}{(rho,u)} in Eulerian coordinate}\label{subsec:holder}
Thus far, we have implicitly used the H\"older continuity of approximate solution: $$|\rho_m|_{1+\sigma,\, 1+\sigma/2}^T,\qu |u_m|_{2+\sigma,\,1+\sigma/2}^T\le C_m.$$ This is used to ensure the validity of translation between Lagrangian coordinate and Eulerian coordinate. Here, the bounding constant for the H\"older norm, $C_m>0$ possibly depends on $m\in\mathbb{N}$. Hence, with the aim of obtaining a solution in the unbounded Eulerian domain $r\in[1,\infty)$, we obtain a H\"older estimates of $(v,u)$ uniformly in $m\in\mathbb{N}$. Moreover, the classical results in \cite{Tani} show that the local-in-time solution to the Navier-stokes equation is H\"{o}lder continuous if so does the initial data. Consequently, the a-priori uniform estimate in the H\"{o}lder norm, which is the main lemma of this subsection, will imply that the limit solution in unbounded Eulerian domain is also H\"{o}lder continuous globally in time, via the continuity argument. 

For simplicity we will abbreviated the parameter $m\in\mathbb{N}$ to denote $(\rho,u)\equiv (\rho_m,u_m)$. In addition, to avoid confusion, we make difference between a solution in the Lagrangian coordinate $(x,t)$ %
and that in the Eulerian coordinate $(r,t)$ %
by writing the former with the hat ``\,$\hat{\ }$\,'' as %
$(\hat{v}, \hat{u}) = (\hat{v}, \hat{u}) (x,t)$ (see Section \ref{sec:RL}).

The uniform H\"older estimate is obtained from the results in Section \ref{sec:higher}, which are restricted in smaller domain of $\mL(T)$. Thus its derivation must also be restricted in its corresponding subdomain translated back to the Eulerian coordinate. For this purpose, we define the sets:
\begin{equation}\label{subdomain}
\Omega_m(T)\vcentcolon=\{ (r,t)\in[1,\infty)\times[0,T] \;\vert\; \eta(X(r,t),t) = 1 \},
\end{equation}
where $X(r,t)$ is defined in \eqref{RInv}, and $\eta$ is the cut-off function in \eqref{eta}. It follows that
\begin{equation*}
\text{if } \ (r,t)\in\Omega_m(T) \ \text{ then } \ X(r,t) \in \supp(\eta).
\end{equation*}
The proof of the next
lemma has already been shown in \cite{NNY} by applying the Schauder theory for the parabolic equations. The main difference for our case is that the domain is restricted in $\Omega_m(T)$ rather than the whole domain $(r,t)\in[1,m]\times[0,T]$. Since the proof is exactly the same, we abbreviate it here. 

\begin{lemma}\label{lemma:holder}
	\label{hth}
	Suppose that the initial data $(\rho_m^0,u_m^0)(r)$ satisfies the H\"older continuity \eqref{He}, for each $m\in\mathbb{N}$. If there exists $T>0$ such that $(\rho, u)(r,t)$ is a solution to \eqref{approx} in the function space \eqref{msp} with the existence time $T$, then there exists a positive constant $C_0>0$ depending
	only the initial data such that for each subset $E\subseteq \Omega_m(T)$,
	\begin{equation*}
		|\rho|_{1+\sigma,\, 1+\sigma/2}^{E},\quad
		|u|_{2+\sigma,\, 1+\sigma/2}^{E} \le C_0 e^{C_0 T}.
	\end{equation*}
\end{lemma}
For other similar proofs, readers are referred to \cite{k-n-z03,k-n81,NNY,NN}. See also \cite{friedman64,lady} for the general theory on this subject.

\section{Solution in the unbounded Eulerian domain}\label{sec:limit}\setcounter{equation}{0}
In this subsection, we take limit $m\to\infty$ to obtain solution $(\rho,u)$ described in Theorem \ref{thm:main}. Since the boundaries of domain in Lagrangian coordinate moves with time, it is best to take this limit $m\to\infty$ in Eulerian coordinate instead. For this purpose, the uniform a-priori estimates derived under Lagrangian coordinate in Sections \ref{sec:apriori}--\ref{sec:higher} are translated back to Eulerian coordinate.

Let $(\rho_m,u_m)(r,t)$ be the approximate solution defined in $(r,t)\in[1,m]\times[0,T]$, given by Lemma \ref{lemma:local}. Since the coordinate transformation described in Section \ref{sec:RL} is constructed using $(\rho_m,u_m)$, the following list of terms appearing in Section \ref{sec:RL} are denoted with subscript $m$ in order to indicate their implicit dependence on the approximation parameter:
\begin{alignat*}{3}
&B_m(t) \equiv B(t), && M_m(t)\equiv M(t), && \text{ for } \ \text{\eqref{B}, \eqref{M},} \\
&R_m(x,t)\equiv R(x,t), \qu && X_m(r,t) \equiv X(r,t), \qu && \text{ for } \ \text{\eqref{RInv}, \eqref{R},}\\
&\mL_m(T)\equiv \mL(T), && S_m(t)\equiv S(t), && \text{ for } \ \text{\eqref{Ldomain},}\\
&M_{m}^0 \equiv M_0, &&  \eta_m(x,t) \equiv \eta(x,t), && \text{ for } \ \text{\eqref{MODE}, \eqref{eta}.} 
\end{alignat*}
For each $m\in\mathbb{N}$, we define $\chi_m\vcentcolon[1,\infty)\times[0,T]\to [0,1]$ as: 
\begin{equation}\label{chi}
\chi_m(r,t) \vcentcolon= %\eta(X(r,t),t) = 
\tilde{\eta}(r-R_m(M_m^0,t)), \qu \text{for } \ (r,t)\in[1,\infty)\times[0,T].
\end{equation}
where $\tilde{\eta}$ is given in \eqref{teta}. Note that $M_m^0\to\infty$ as $m\to\infty$. It follows that
\begin{proposition}\label{prop:chi}
For each $m\in\mathbb{N}$, if $(x,t)\in \mL_m(T)$ then $\chi_m(R_m(x,t),t)=\eta_m(x,t)$. In addition, for each $(r,t)\in[1,\infty)\times[0,T]$,
\begin{equation*}
\lim\limits_{m\to\infty}(\chi_m,\partial_r\chi_m)(r,t) = (1,0),\qu 0\le \chi_m\le 1, \qu \lp{\infty}{\partial_r\chi_m}\le \sqrt{8}.
\end{equation*}
Moreover, for each $f\in L^\infty(0,T;L^1(1,\infty))$,
\begin{align*}
\lim\limits_{m\to\infty}\int_{0}^{T}\!\!\!\int_{1}^{\infty} |\partial_t \chi_m|^2 f \, \dif r \dif t = 0 \qu \text{and} \qu  \int_0^{T} \sup\limits_{r\ge 1} | \partial_t \chi_m(r,t) |^2 \dif t \le C_0(1+T),
\end{align*}
where $C_0>0$ is a constant independent of $m\in\mathbb{N}$ and $T>0$.
\end{proposition}
\begin{proof}
First, we show that $\inf_{t\in[0,T]}R_m(M_{m}^0,t)\to \infty$. Using (\ref{B}), (\ref{R}), and Lemma \ref{lemma:v}, it follows that as $m\to\infty$,
\begin{align}\label{temp:chi1}
&\inf_{t\in[0,T]}R_m(M_{m}^0,t)= \inf_{t\in[0,T]} \Big( 1 + n\int_{B_m(t)}^{M_{m}^0}v_m(x,t) \dif x \Big)^{\frac{1}{n}}\nonumber\\ &\ge \inf_{t\in[0,T]} \big( 1 + n \ul{v} ( M_{m}^0 - B_m(t) ) \big)^{\frac{1}{n}} \ge \big( 1 + n \ul{v} ( M_{m}^0 - \ol{v} |u_b| T ) \big)^{\frac{1}{n}} \to \infty
\end{align}
Due to (\ref{temp:chi1}), for a given $(r,t)\in [1,\infty)\times[0,T]$, there exists $N=N(r,t)\in \mathbb{N}$ such that if $m\ge N$ then $r\le R_m(M_m^0,t)-1$. Thus by definition (\ref{chi}), $\chi_m(r,t)=1$ for all $m\ge N$, which implies that $\chi_m(r,t) \to 1 $ as $m\to\infty$ for all $(r,t)\in [1,\infty)\times[0,T]$.

Next, taking spatial derivative on $\chi_m(r,t)$, we have $\partial_r \chi_m(r,t) = \tilde{\eta}^{\prime}(r-R_m(M_m^0,t))$. From the construction (\ref{teta}), $\supp(\tilde{\eta}^{\prime})\subseteq [-1,0]$, hence combining with the limit (\ref{temp:chi1}), we have that for all $(r,t)\in[1,\infty)\times[0,T]$,
\begin{equation}\label{temp:chi2}
\partial_r \chi_m(r,t) = \tilde{\eta}^{\prime}(r-R_m(M_m^0,t)) \to 0 \ \text{ as } \ m\to\infty.
\end{equation}

Taking temporal derivative on $\chi_m(r,t)$ it follows from (\ref{Rdiff}) that, 
\begin{equation}\label{temp:chi3}
\partial_t \chi_m(r,t) = - u_m(R_m(M_m^0,t),t) \cdot \tilde{\eta}^{\prime}(r-R_m(M_m^0,t)).
\end{equation}
Denoting $ \hat{u}_m(x,t) \equiv u_m(R_m(x,t),t)$. Then by (\ref{teta}) and Corollary \ref{corol:infpsi}, we have
\begin{equation*}
\int_{0}^{T}\sup_{r\ge 1}|\partial_t \chi_m(r,t)|^2 \dif t \le 8 \int_{0}^{T} |\hat{u}_m(M_m^0,t)|^2 \sup_{r\ge 1} \chi_m(r,t) \dif t \le C_0(1+T).
\end{equation*}
Set $f\in L^{\infty}(0,T;L^1(1,\infty))$. Using (\ref{temp:chi2})--(\ref{temp:chi3}), Corollary \ref{corol:infpsi}, and
the
dominated convergence theorem we obtain that
\begin{align*}
\Big| \int_{0}^{T}\!\!\!\int_{1}^{\infty} f |\partial_t \chi_m|^2 \dif r \dif t \Big| \le& \Big( \int_{0}^{T}\!\! |\hat{u}_m(M_m^0,t)|^2\, \dif t \Big)  \Big( \sup_{t\in[0,T]} \int_{1}^{\infty}\! |\partial_r\chi_m|^2 |f| \, \dif r \Big)\\
\le & C_0(1+T) \Big( \sup_{t\in[0,T]} \int_{1}^{\infty}\! |\partial_r\chi_m|^2 |f| \, \dif r \Big) \to 0 \ \ \text{as } \ m\to\infty.
\end{align*}
This concludes the proof.
\end{proof}
Using the coordinate transformation \eqref{rR}, we translate the uniform a-priori estimates obtained in Sections \ref{sec:approx}--\ref{sec:higher} to the bounded Eulerian coordinate $(r,t)\in[1,m]\times[0,T]$. As a result, we have the following proposition:
\begin{proposition}\label{prop:mEulEst}
There exists constants $0<\ul{\rho}<\bar{\rho}<\infty$ which only depends on initial data such that
\begin{equation*}
	\ul{\rho} \le \rho_m(r,t) \le \bar{\rho}, \qu \text{for all } \  (r,t)\in [1,m]\times[0,T].
\end{equation*}
Moreover, there exists a positive constant $C>0$ independent of $m\in\mathbb{N}$ and $T>0$ such that for all $m\in\mathbb{N}$,
\begin{align*}
&\sup\limits_{t\in[0,T]} \int_{1}^{m} \{ (\rho_m-\rt)^2 + (\bar{u}_m-\ut)^2  \}(r,t) r^{n-1} \dif r  \le C E_0,\\
&|u_b|\int_{0}^{T} \big\{ |\rho_m(1,t)-\rt(1)|^2 + (\rho_m-\rt)_r^2(1,t) \big\} \, \dif t \le C E_0\\
&\int_{0}^{T}\!\!\!\int_{1}^{m}\!\! \Big\{ (u_m-\ut)_r^2 + \dfrac{(u_m-\ut)^2}{r^2} + |u_b| \dfrac{(u_m-\ut)^2}{r^n} + |u_b|^3 \dfrac{(\rho_m-\rt)^2}{r^{3n-2}} \Big\} r^{n-1} \dif r \dif t \le C E_0,\\
&\sup\limits_{t\in[0,T]}\int_{1}^{m} \chi_m |(\rho_m-\rt)_r|^2(r,t) r^{n-1} \dif r + \int_{0}^{T}\!\!\int_{1}^{m}\chi_m |(\rho_m-\rt)_r|^2 r^{n-3} \dif r \dif t\le C E_0,\\
&\sup\limits_{t\in[0,T]}\int_{1}^{m} \chi_m |(u_m-\ut)_r|^2(r,t) r^{n-1} \dif r + \int_{0}^{T}\!\!\int_{1}^{m}\chi_m |(u_m-\ut)_{rr}|^2 r^{n-3} \dif r \dif t\le C E_0,\\
&\int_{0}^{T}\!\!\int_{1}^{m} |\partial_t\rho_m|^2 r^{n-1} \dif r \dif t \le C E_0, \qu \int_{0}^{T}\!\!\int_{1}^{m} \chi_m |\partial_t u_m|^2 r^{n-1} \dif r \dif t \le CE_0(1+T).
\end{align*}
\end{proposition}

Next, the solution $(\rho_m,u_m)$ is extended to the whole domain $(r,t)\in[1,\infty)\times[0,T]$ as follows: let $\varphi_m$ be the cut-off function introduced in \eqref{cut1}. For all $(r,t)\in[1,\infty)\times[0,T]$, we set 
\begin{equation}\label{ext}
	\breve{\rho}_m\vcentcolon= \rho_m \varphi_m + \rt(1-\varphi_m), \qu \breve{u}_m\vcentcolon= u_m \varphi_m + \ut (1-\varphi_m).
\end{equation}
From Proposition \ref{prop:mEulEst} and \eqref{ext}, one can verify that  
\begin{equation*}
\ul{\rho} \le \breve{\rho}_m(r,t) \le \bar{\rho}, \qu \text{for all } \  (m,r,t)\in\mathbb{N}\times[1,\infty)\times[0,T],
\end{equation*}
for some constants $0<\ul{\rho}<\bar{\rho}<\infty$ which only depends on initial data. Moreover, there exists a positive constant $C>0$ independent of $m\in\mathbb{N}$ and $T>0$ such that for all $m\in\mathbb{N}$,
\begin{align*}
&\sup\limits_{t\in[0,T]} \int_{1}^{\infty} \{ (\breve{\rho}_m-\rt)^2 + (\breve{u}_m-\ut)^2  \}(r,t) r^{n-1} \dif r \le C E_0,\\
&\sup\limits_{t\in[0,T]}\int_{1}^{\infty} \chi_m \{  |(\breve{\rho}_m-\rt)_r|^2 + |(\breve{u}_m-\ut)_r|^2\}(r,t) r^{n-1} \dif r \le C E_0,\\
&\int_{0}^{T}\!\!\int_{1}^{\infty} |\partial_t \breve{\rho}_m|^2 r^{n-1} \dif r \dif t \le C E_0, \qu \int_{0}^{T}\!\!\int_{1}^{\infty} \chi_m |\partial_t \breve{u}_m|^2 r^{n-1} \dif r \dif t \le CE_0(1+T),\\
&\int_{0}^{T}\!\!\int_{1}^{\infty} \chi_m|(\breve{u}_m-\ut)_{rr}|^2 r^{n-3} \dif r \dif t\le C E_0.
\end{align*}
By the weak and weak-$\ast$ compactness of
the Sobolev spaces and Proposition \ref{prop:chi}, there exists some function $(\rho,u)\in H^1_{\text{loc}}$ such that as $m\to\infty$,
\begin{align*}
	&r^{\frac{n-1}{2}}\chi_m(\breve{\rho}_m - \rt, \breve{u}_m-\ut) \overset{\ast}{\rightharpoonup} r^{\frac{n-1}{2}}(\rho - \rt, u-\ut) && \text{(weakly-$*$) } && \text{in } \ L^\infty(0,T; H^1(0, \infty)),
	\\
	&r^{\frac{n-1}{2}}\chi_m( \partial_t\breve{\rho}_{m}, \partial_t\breve{u}_{m}) \rightharpoonup r^{\frac{n-1}{2}}(\rho_t, u_t) && \text{(weakly) } && \text{in } \ L^2 (0,T; L^2(0, \infty)),
	\\
	&r^{\frac{n-3}{2}}\chi_m (\breve{u}_{m}-\ut)_{rr} \rightharpoonup r^{\frac{n-3}{2}} (u-\ut)_{rr} && \text{(weakly) } && \text{in } \ L^2 (0,T; L^2(0, \infty)).
\end{align*}
Combining the above convergence with uniform H\"older estimate Lemma \ref{lemma:holder}, one can show by
the
Arzel\`a-Ascoli theorem that there exists a realisation of the limit function $(\rho,u)$ such that 
\begin{equation*}
|\rho|_{1+\sigma,\, 1+\sigma/2}^{Q_T},\quad |u|_{2+\sigma,\, 1+\sigma/2}^{Q_T} \le C_0 e^{C_0 T},
\end{equation*}
where $C_0>0$ is a constant depending only on the initial data, and $Q_T\vcentcolon=[1,\infty)\times[0,T]$. Furthermore, we see that $(\rho,u)$ is a unique solution in $(r,t)\in [1,\infty) \times [0,T]$
to the problem \eqref{nse}--\eqref{compa}.
%\begin{subequations}\label{xx1}
%\begin{gather}
%\rho-\rt, \; u - \ut , \; r^{\frac{n-1}{2}} (\rho-\rt)_r, \; r^{\frac{n-1}{2}} (u-\ut)_r \in C ([0,T]; L^2(0, \infty)), \\
%r^{-1} u, \; r^{\frac{n-3}{2}} (\rho-\rt)_r, \; r^{\frac{n-3}{2}} (u-\ut)_{rr} \in L^2(0,T; L^2(0, \infty)),
%\end{gather}
%\end{subequations}
Moreover, utilizing Proposition \ref{prop:mEulEst}, we obtain
\begin{subequations}\label{finalEst}
\begin{gather}
\sup_{t\in[0,T]}\int_1^\infty \{ (\rho-\rt)^2 + (u-\ut)^2 + (\rho-\rt)_r^2 + (u-\ut)_r^2 \}(r,t) r^{n-1} \,  \dif r \le C E_0, \label{finalEst1} \\
|u_b|\int_{0}^{T} \big\{ |\rho(1,t)-\rt(1)|^2 + (\rho-\rt)_r^2(1,t) \big\} \, \dif t \le C E_0,\\
\int_0^T \!\!\! \int_1^\infty \big\{\dfrac{(u-\ut)^2}{r^2} + |u_b|^3\dfrac{(\rho-\rt)^2}{r^{3n-2}} \big\} r^{n-1}  \, \dif r \dif t
\le C E_0,
\label{finalEst2}\\
\int_0^T \!\!\! \int_1^\infty \big\{ (u-\ut)_r^2 + \dfrac{(\rho-\rt)_r^2}{r^2} + \dfrac{(u-\ut)_{rr}^2}{r^2} + u_t^2 + \rho_t^2 \big\} r^{n-1}  \, \dif r \dif t
\le C E_0,\label{finalEst3}
\end{gather}
\end{subequations}
where
$C$ is a positive constant depending only on the initial data \eqref{ic}.

\section{Large time behavior of solutions in Eulerian coordinate}\label{sec:tAsymp}\setcounter{equation}{0}
Using the estimate \eqref{finalEst}, we obtain the asymptotic stability of the stationary solution $(\rt,\ut)$ in the Eulerian coordinate.
Precisely, we show the convergence \eqref{tcvg}.
\begin{proposition}\label{prop:tAsymp}
The convergence \eqref{tcvg} holds under the same assumptions as in Theorem \ref{thm:main}.
\end{proposition}
\begin{proof}
At first, we prove the assertion concerning $u$, that is, %
 $\lp{\infty}{u(\cdot,t)-\tilde{u}(\cdot)} \to 0$ as $t \to \infty$.
Applying the Schwarz inequality and using \eqref{finalEst}, we see that
\begin{align*}
&\sup_{r \in [1,\infty)} (u-\ut)^2\\
\le&
C \Big( \int_1^\infty \frac{(u-\ut)^2}{r^{n-1}} \, \dif r \Big)^{1/2}
\Big( \int_1^\infty r^{n-1} (u-\ut)_r^2 \, \dif r \Big)^{1/2}
\le
C \Big( \int_1^\infty \frac{(u-\ut)^2}{r^{n-1}} \, \dif r \Big)^{1/2}.
\end{align*}
Thus, in order to prove this assertion, it suffices to show
\begin{equation}
I_1(t)
:=
\int_1^\infty \frac{(u-\ut)^2}{r^{n-1}}(r,t) \, \dif r
\to 0
\quad \text{as} \quad t \to \infty.
\label{m2}
\end{equation}
To this end, it suffices to
 show that $\frac{\dif}{\dif t} I_1(t) \in L^1(0,\infty)$
%which immediately gives (\ref{m2})
as  $I_1 (t) \in L^1(0, \infty)$ owing to (\ref{finalEst}).
Differentiating $I_1(t)$ in $t$, we have
\begin{align}
&\Big|\frac{\dif}{\dif t} I_1(t)\Big|
=
2\int_1^\infty |u_t (u-\ut)|  r^{1-n} \dif r \le \int_{1}^{\infty}\!\! |u_t|^2 r^{n-1}\, \dif r + \int_{1}^{\infty}\!\! \dfrac{|u-\ut|^2}{r^2} r^{n-1}\, \dif r.
\label{m3}
\end{align}
By (\ref{finalEst}), the above inequality implies that $\frac{d}{dt} I_1 (t)$ and $I_1(t)$ are both integrable over $t\in(0,\infty)$. Consequently, we obtain the convergence $\lp{\infty}{u(\cdot,t)-\ut(\cdot)} \to 0$ as $t\to\infty$.

Next, we show $\lp{\infty}{\rho(\cdot,t) - \rt(\cdot)} \to 0$ as $t \to \infty$. By the similar argument to the above, we have
\begin{align*}
&\sup_{r \in [1,\infty)} |\rho(r,t)-\rt(r)|^2\\
\le& |\rho(1,t) - \rt(1) |^2  +
2\Big( \int_1^\infty\!\! (\rho-\rt)^2r^{n-1} \, \dif r \Big)^{1/2}
\Big( \int_1^\infty\!\! \dfrac{(\rho-\rt)_r^2}{r^{n-1}} \, \dif r \Big)^{1/2}\\
\le& |\rho(1,t) - \rt(1) |^2 +
C E_0 \Big( \int_1^\infty\!\! \dfrac{(\rho-\rt)_r^2}{r^{n-1}} \, \dif r \Big)^{1/2}.
\end{align*}
It suffices to show that as $t\to\infty$,
\begin{equation}
J_1(t)\vcentcolon=|\rho(1,t) - \rt(1) |^2 \to 0 ,\qquad J_2(t)
\vcentcolon=
\int_1^\infty\!\! \dfrac{(\rho-\rt)_r^2}{r^{n-1}} (r,t)  \, \dif r 
\to 0.
\label{m12}
\end{equation}
First, we consider the term $J_1(t)$. From (\ref{finalEst}), we have $J_1(t)\in L^{1}(0,\infty)$. In addition, taking $r\to 1$ on the equation (\ref{nse1}), we have
\begin{multline*}
\rho_t (1,t) = - u_b  \big(\rho_r(1,t)-\rt_r(1)\big) - \rho(1,t) \big(u_r(1,t)-\ut_r(1)\big)\\
- \ut_r(1) \big(\rho(1,t)-\rt(1)\big) - (n-1) u_b \big(\rho(1,t)-\rt(1)\big). 
\end{multline*}
Using this, and
the Sobolev embedding theorem, it follows that
\begin{align*}
&\Big|\dfrac{\dif }{\dif t} J_1(t)\Big| = 2\big| \big(\rho(1,t)-\rt(1)\big) \d_t\rho(1,t) \big|\\
\le & C |\rho(1,t)-\rt(1)|^2 + C |\rho_r(1,t)-\rt_r(1)|^2 + C\int_{1}^{\infty}\!\! \Big\{ \dfrac{(u-\ut)_{rr}^2}{r^2} + (u-\ut)_r^2 \Big\} r^{n-1}\, \dif r 
\end{align*}
By (\ref{finalEst}) and the above inequality, $J_1(t)$ and $\frac{\dif}{\dif t} J_1(t)$ are both in $L^1(0,\infty)$, which implies that $J_1(t)\to 0$ as $t\to\infty$. Next we consider $J_2(t)$. By a straightforward computation on (\ref{nse1}), we see that
%\begin{align*}
%\rho_{t}
%=& \rho \big\{ (u-\ut)_r + \dfrac{n-1}{r}(u-\ut) \big\} + (\rho-\rt)_r(u-\ut) + \rt_r(u-\ut)\\
%&+ \dfrac{\d_r (r^{n-1}\ut)}{r^{n-1}} (\rho-\rt) + \ut (\rho-\rt)_r.
%\label{m6}
%\end{align*}
\begin{align*}
\rho_{tr} =& u (\rho-\rt)_{rr}  + \rho (u-\ut)_{rr} + 2 (u-\ut)_r(\rho-\rt)_r + \dfrac{n-1}{r} (u-\ut) (\rho-\rt)_r\\
&+ \big\{ 2 \d_r \ut + \dfrac{n-1}{r} \ut \big\} (\rho-\rt)_r + \big\{ 2\rt_r + \dfrac{n-1}{r}\rho \big\} (u-\ut)_r - \dfrac{n-1}{r^2}(\rho-\rt)(u-\ut) \\
&+\big\{ \rt_{rr} + \dfrac{n-1}{r}\rt_r - \dfrac{n-1}{r^2} \rt \big\} (u-\ut) + \d_r\big( \dfrac{\d_r(r^{n-1}\ut)}{r^{n-1}} \big) (\rho-\rt).
\end{align*}
By the above identity, integration by parts, and 
the
Sobolev embedding theorem, we have:
\begin{align*}
&\Big|\dfrac{\dif}{\dif t} J_2(t) \Big|
=
2 \Big| \int_1^\infty\!\! \rho_{tr} (\rho-\rt)_r  r^{1-n} \, \dif r \Big|\\
\le &C \int_{1}^{\infty}\!\! \Big\{ \dfrac{(u-\ut)_{rr}^2}{r^2} +  (u-\ut)_r^2 + (\rho-\rt)_r^2  + \dfrac{|u-\ut|^2}{r^2} + \dfrac{|\rho-\rt|^2}{r^{3n-2}} \Big\} r^{n-1} \,\dif r. 
\end{align*}
Thus $\frac{d}{dt} J_2 (t)$ and $J_2(t)$ are both integrable over $t\in(0,\infty)$. Hence $J_2(t)\to 0$ as $t\to\infty$, which concludes that $\lp{\infty}{\rho(\cdot,t) - \rt(\cdot)} \to 0$ as $t \to \infty$.
\end{proof}

\appendix

\section{Proof for Proposition \ref{prop:phiG} }\label{appen:G}
In this appendix, we give the proof for Proposition \ref{prop:phiG}. 
By the fact that $p(v)=Kv^{-\gamma}$, we have that
\begin{equation*}
    G(v,w) = \int_{w^{-1}}^{v^{-1}} \dfrac{p(z)-p(w)}{z^2} \dif z = \left\{
    \begin{aligned}
    &\dfrac{K}{\gamma-1} (v^{1-\gamma}-w^{1-\gamma}) + K w^{-\gamma} (v-w)  && \text{if} \qu \gamma>1,\\
    &K\big( \dfrac{v}{w} - 1 -\log \dfrac{v}{w} \big) && \text{if} \qu \gamma =1.
    \end{aligned}\right.
\end{equation*}
It is verified that for $\gamma\ge 1$, $G(y,y)=\partial_v G(y,y) = \partial_w G(y,y) =0$ for all $y > 0$. %Denote $X=(v,w)$, $Y=(w,w)$, and $\nabla G \equiv (\partial_v G, \partial_w G )^{\top}$. 
By Taylor's theorem, there exists $\theta=\theta(v,w)\in(0,1)$ such that if one sets the vector $\xi\vcentcolon = (\theta v + (1-\theta)w , w)$, then it holds
\begin{align*}
G(v,w) %= G(Y) + (X-Y) \cdot \nabla G(Y) + \dfrac{1}{2} (X-Y) \cdot \begin{pmatrix} \partial_v^2 G(\xi) & \partial_v\partial_w G(\xi)\\ \partial_w\partial_v G(\xi) & \partial_w^2 G(\xi)  \end{pmatrix} \cdot (X-Y)^{\top}\\
= \dfrac{1}{2} (v-w,0) \cdot \begin{pmatrix} \partial_v^2 G(\xi) & \partial_v\partial_w G(\xi)\\ \partial_w\partial_v G(\xi) & \partial_w^2 G(\xi)  \end{pmatrix}
 \cdot 
\begin{pmatrix} v-w \\ 0 \end{pmatrix}.
\end{align*}
Thus computing the Hessian of $G$, we obtain that for $\gamma \ge 1$
\begin{align*}
G(v,w) = %\left\{ \begin{aligned} &
\dfrac{K\gamma}{2}  \dfrac{|v-w|^2}{(\theta v + (1-\theta) w )^{1+\gamma}} %&& \text{if} \qu \gamma>1,\\
%&\dfrac{K}{2} \dfrac{|v-w|^2}{(\theta v + (1-\theta)w)^2} && \text{if} \qu \gamma=1.
%\end{aligned} \right.
\qu \text{for  $v>0$ and $w\ge 0$}. 
\end{align*}
If $v\le \tilde{v}$, then setting $w=\tilde{v}$ in the above we obtain that $G(v,\tilde{v}) \ge (K\gamma/2) \vt^{-\gamma-1} |v-\vt|^2$. This proves the first inequality in Proposition \ref{prop:phiG}. 

Next we define the function $H(\rho,\sigma)\vcentcolon= \rho G(\rho^{-1},\sigma^{-1})$, then it can also be verified that $H(y,y)=\partial_\rho H(y,y) = \partial_\sigma H(y,y) = 0$ for $y\ge 0$. By the same argument as before, there exists $\lambda=\lambda(\rho,\sigma) \in (0,1)$ such that if one sets $\zeta \vcentcolon= (\lambda\rho + (1-\lambda)\sigma,\sigma)$ then
\begin{align*}
H(\rho,\sigma)
= \dfrac{1}{2} (\rho-\sigma,0) \cdot \begin{pmatrix} \partial_\rho^2 H(\zeta) & \partial_\sigma\partial_\rho H(\zeta)\\ \partial_\rho\partial_\sigma H(\zeta) & \partial_\sigma^2 H(\zeta)  \end{pmatrix} \cdot
\begin{pmatrix}
\rho-\sigma \\ 0 \end{pmatrix}.
\end{align*}
Computing the Hessian of $H$, we obtain that for $\gamma \ge 1$,
\begin{equation*}
    H(\rho,\sigma) = %\left\{ \begin{aligned}&
    \dfrac{K \gamma}{2} (\lambda \rho + (1-\lambda)\sigma )^{\gamma-2} |\rho-\sigma|^2 %&& \text{if} \qu \gamma>1,\\
    %&\dfrac{K}{2} \dfrac{|\rho-\sigma|^2}{\lambda \rho + (1-\lambda)\sigma} && \text{if} \qu \gamma=1,
    %\end{aligned}\right.
\end{equation*}
For the case $v>\vt$, we have $\rho = 1/v < 1/\vt =\rt$. Thus setting $\sigma=\rt$ in the above equation, it follows that if $1\le \gamma\le 2$ then $H(\rho,\rt) \ge (K\gamma/2) \rt^{\gamma-2} |\rho-\rt|^2 = (K\gamma/2) \vt^{-\gamma} v^{-2} |v-\vt|^2 $. Hence we have $G(v,\vt) \ge (K\gamma/2) \vt^{-\gamma} v^{-1} |v-\vt|^2$, which is the second inequality of Proposition \ref{prop:phiG}. \null\hfill\qedsymbol

%\bibliographystyle{siam}
%\nocite{*}
%\bibliography{refs}

\end{document}